%
%
%
%
%

\input amstex  

\input amssym
\input amssym.def

\magnification 1200
\loadmsbm
\parindent 0 cm


\define\nl{\bigskip\item{}}
\define\snl{\smallskip\item{}}
\define\inspr #1{\parindent=20pt\bigskip\bf\item{#1}}
\define\iinspr #1{\parindent=27pt\bigskip\bf\item{#1}}
\define\ainspr #1{\parindent=24pt\bigskip\bf\item{#1}}
\define\aiinspr #1{\parindent=31pt\bigskip\bf\item{#1}}

\define\einspr{\parindent=0cm\bigskip}

\define\ot{\otimes}

\centerline{\bf Multiplier Hopf algebroids arising from weak multiplier Hopf algebras}
\bigskip
\centerline{\it T.\ Timmermann \rm $^{(1)}$ and \it A.\ Van Daele \rm $^{(2)}$} 
\nl\nl\nl
{\bf Abstract} 
\nl 
It is well-known that any weak Hopf algebra gives rise to a Hopf algebroid. Moreover it is possible to characterize those Hopf algebroids that arise in this way.
\snl
Recently, the notion of a weak Hopf algebra has been extended to the case of algebras without identity. This led to the theory of weak multiplier Hopf algebras. Similarly also the theory of Hopf algebroids was recently developed for algebras without identity. They are called multiplier Hopf algebroids. Then it is quite natural to investigate the expected link between weak multiplier Hopf algebras and multiplier Hopf algebroids. This relation has been considered already in the original paper on multiplier Hopf algebroids. In this note, we investigate the connection further. 
\snl
First we show that any regular weak multiplier Hopf algebra gives rise, in a natural way, to a regular multiplier Hopf algebroid. Secondly we give a characterization, mainly in terms of the base algebra, for a regular multiplier Hopf algebroid to have an underlying weak multiplier Hopf algebra. We illustrate this result with some examples. In particular, we give examples of multiplier Hopf algebroids that do not arise from a weak multiplier Hopf algebra.
\nl\nl
{\it June 2014} (Version 1.1)
\vskip 5 cm
\hrule
\bigskip
\parindent 0.7 cm
\item{($1$)} FB Mathematik und Informatik, University of M\"unster, Einsteinstrasse 62, D-48149 M\"unster, Germany. {\it E-mail}: timmermt\@math.uni-muenster.de
\item{($2$)} Department of Mathematics, University of Leuven, Celestijnenlaan 200B,
B-3001 Heverlee, Belgium. {\it E-mail}: Alfons.VanDaele\@wis.kuleuven.be

\parindent 0 cm

\newpage

\bf 0. Introduction \rm
\nl
Let $G$ be a {\it finite groupoid}. Consider the algebra $K(G)$ of all complex functions on $G$ with point-wise operations. It will be denoted by $A$. Define a coproduct $\Delta:A\to A\ot A$ by
$$\Delta(f)(p,q)= 
	\cases
		f(pq) & \text{if } pq \text{ is defined},\\
				0 & \text{otherwise}
	\endcases \tag"(0.1)"$$
where $f\in A$ and $p,q\in G$. We identify $A\ot A$ with the algebra $K(G\times G)$ of all complex functions on the Cartesian product $G\times G$ of $G$ with itself. It is not hard to see that $\Delta$ is a homomorphism. It is also coassociative in the sense that $(\Delta\ot\iota)\Delta=(\iota\ot\Delta)\Delta$ where $\iota$ denotes the identity map. The algebra $A$ clearly has an identity. However, the coproduct $\Delta$ need not be unital. This will be the case only when the product $pq$ is defined for all pairs $(p,q)\in G\times G$, i.e.\ when $G$ is a group. The pair $(A,\Delta)$ is a {\it weak Hopf algebra} as introduced in [B-N-S].
\snl
If $G$ is {\it no longer} assumed to be {\it finite}, we take for $A$ the algebra $K(G)$ of complex functions with {\it finite support} instead of all complex functions. This is now an algebra without identity. The product is non-degenerate and we can consider the multiplier algebra $M(A)$. It is naturally identified with the algebra $C(G)$ of all complex functions on $G$. We can still define a coproduct $\Delta$ as in the formula (0.1) above but now $\Delta$ maps $A$ to $M(A\ot A)$, the multiplier algebra of $A\ot A$. Again $A\ot A$ is identified with $K(G\times G)$, the algebra of complex functions with finite support in $G\times G$, while $M(A\ot A)$ is identified with $C(G\times G)$, the algebra of all complex functions on $G\times G$. In this case, the pair $(A,\Delta)$ is a {\it weak multiplier Hopf algebra} as introduced and studied in [VD-W1].
\nl
Next consider any {\it weak Hopf algebra} $(A,\Delta)$. It has a unique antipode $S$. Denote by $\varepsilon_s$ and $\varepsilon_t$ respectively the source and target maps, defined from $A$ to itself by
$$\varepsilon_s(a)=\sum S(a_{(1)}) a_{(2)}
\quad\quad\text{and}\quad\quad
\varepsilon_t(a)=\sum a_{(1)}S(a_{(2)})$$
for $a\in A$. We use the Sweedler notation $\sum a_{(1)}\ot a_{(2)}$ for $\Delta(a)$. The ranges of $\varepsilon_s$ and $\varepsilon_t$ will be denoted by $B$ and $C$ respectively. Then $B$ and $C$ are two commuting subalgebras of $A$. If we use $S_B$ and $S_C$ for the restrictions of the antipode $S$ to these subalgebras $B$ and $C$ respectively, it turns out that $S_B$ is an anti-isomorphism from $B$ to $C$ and that $S_C$ is an anti-isomorphism from $C$ to $B$. Observe that in general, they are not inverses of each other. 
\snl
By assumption, here the algebra $A$ is an algebra with an identity $1$. We have that $\Delta(1)\in B\ot C$ and that $B$ and $C$ are the smallest subspaces of $A$ with this property. In other words, $B$ is the left leg of $\Delta(1)$ and $C$ is the right leg of $\Delta(1)$. Recall that $\Delta(1)$ is not necessarily equal to $1\ot 1$. If that is the case, that is when $(A,\Delta)$ is actually a Hopf algebra, we have that $B$ and $C$ are just the scalar multiples of the identity and the source and target maps are essentially the counit of the Hopf algebra $(A,\Delta)$.
\snl
Consider now two {\it balanced tensor products}. The first one, denoted as $A\ot_\ell A$, is characterized by the requirement that 
$xa\ot b=a\ot S_B(x)b$ for all $a,b\in A$ and $x\in B$. In the second one, denoted as $A\ot_r A$, we have $a\ot by=aS_C(y)\ot b$ for all $a,b\in A$ and $y\in C$.
Denote by $\pi_\ell$ and $\pi_r$ the canonical quotient maps from $A\ot A$ to $A\ot_\ell A$ and from $A\ot A$  to $A\ot_r A$ respectively. Define
$$\Delta_B(a)=\pi_\ell (\Delta(a))
\qquad\quad\text{and}\qquad\quad
\Delta_C(a)=\pi_r(\Delta(a))$$
for $a\in A$. The maps 
$$\Delta_B:A\to A\ot_\ell A
\qquad\quad\text{and}\qquad\quad
\Delta_C:A\to A\ot_r A$$
are a left and and a right coproduct that turn the above data $(B,C,A,S_B,S_C)$ into a {\it Hopf algebroid}.
\snl
There is much more to explain about the construction above. We will do this further in the paper where in fact, we describe the construction in the more general case of a regular weak multiplier Hopf algebra. 
\nl
This is the main topic of this work and it takes us to the content of the paper.
\nl
\it Content of the paper \rm
\nl
{\it Section} 1 is a preliminary section. It mainly contains  a review of the notion of a {\it regular weak multiplier Hopf algebra}, one of the basic ingredients of this paper and some of its basic properties. Special attention is given to the source and target maps and the source and target algebras. We do not recall the notion of a {\it multiplier Hopf algebroid} in this preliminary section. As it turns out, it will be easier to do it {\it on the way}, in Section 3, where we pass from a regular weak multiplier Hopf algebra to a multiplier Hopf algebroid. 
\snl
{\it Section} 2 is devoted to what is called in [T-VD] {\it quantum graphs} and the relation with separability idempotents as studied in [VD2]. This is the basic material, commonly needed in both Section 3 and Section 4. If a Hopf algebroid is coming from a weak Hopf algebra as explained above, there is a left quantum graph $(B,A,\iota_B,S_B)$, together with a right quantum graph $(C,A,\iota_C,S_C)$ and the two are compatible. The separability idempotent in this case will be nothing else but $\Delta(1)$. In this section we consider the generalizations of these objects to the framework of weak multiplier Hopf algebras. We will not yet consider the coproducts. That will only be done in the next two sections. 
\snl
In {\it Section} 3 we show how we can associate a multiplier Hopf algebroid to any regular weak multiplier Hopf algebra. The procedure is quite standard and there are no surprises. It is indeed very similar to the same construction in the case of a weak Hopf algebra, with the typical problems due to the fact that the algebras have no identity and that the coproducts do not map into the tensor product, but into a bigger space.
During the process, we recall the various notions that we encounter in the definition of a multiplier Hopf algebroid. We finish the section by looking at the concrete expressions for the counit, as considered in the theory of multiplier Hopf algebroids, in terms of the counital maps, considered in the theory of weak multiplier Hopf algebras. We also illustrate some formulas involving the antipode in the two settings.
\snl
In {\it Section} 4 we consider the opposite direction. We start with a regular multiplier Hopf algebroid and we discuss the necessary and sufficient conditions for it to come from a weak multiplier Hopf algebra as in the previous section. 
\snl
The main condition is that the base algebra $B$ is separable Frobenius. This means that there is a separability idempotent $E\in M(B\ot C)$ in the sense of [VD2]. However, this is not sufficient. The anti-isomorphisms $S_B:B\to C$ and $S_C:C\to B$, that are given for the multiplier Hopf algebroid have to coincide with the anti-isomorphisms induced by $E$. In other words, we have to require that
$$E(x\ot 1)=E(1\ot S_B(x))
\qquad\quad\text{and}\qquad\quad
(1\ot y)E=(S_C(y)\ot 1)E$$
for all $x\in B$ and $y\in C$ (where $S_B$ and $S_C$ are the original anti-isomorphisms that come with the definition of the multiplier Hopf algebroid).
\snl
Still, this will not be sufficient. We also need the equality of the counits $\varepsilon$ and $\varepsilon'$ on $A$ defined by
$$\varepsilon=\varphi_B\circ\varepsilon_B
\qquad\quad\text{and}\qquad\quad
\varepsilon'=\varphi_C\circ\varepsilon_C.$$
Here $\varphi_B$ and $\varphi_C$ are the unique linear functionals on $B$ and $C$ respectively characterized by 
$$(\varphi_B\ot\iota)E=1
\qquad\quad\text{and}\qquad\quad
(\iota\ot\varphi_C)E=1$$
and $\varepsilon_B$ and $\varepsilon_C$ are the counital maps that are obtained within the theory for the original multiplier Hopf algebroid.
\snl
In this section, we further discuss these conditions and compare them with the conditions found for a Hopf algebroid to arise from a weak Hopf algebra as studied in [B]. 
\snl
In {\it Section} 5 we discuss some examples and special cases.
\snl
In the last section, {\it Section} 6, we draw some conclusions and discuss possible future research related with the results obtained in this paper.
\nl
\it Conventions and basic references \rm
\nl
We only consider algebras over the field of complex numbers (although we believe that most of our results are still valid for more general fields). We do not require our algebras to have a unit, but we do need that the product is non-degenerate (as a bilinear form). Our algebras are all idempotent, either by assumption, or by a result. In fact, they even will have local units and this implies that the product is non-degenerate and that the algebra is idempotent.
\snl
For algebras with a non-degenerate product, it is possible to define the multiplier algebra $M(A)$. It can be characterized as the largest algebra with identity containing $A$ as an essential two-sided ideal. We use $1$ for the unit in the multiplier algebras. Of course $M(A)=A$ if and only if $A$ already has an identity. If $A$ is non-degenerate, then so is the tensor product $A\ot A$ and we consider also its multiplier algebra $M(A\ot A)$. We have natural imbeddings
$$A\ot A\subseteq M(A)\ot M(A)\subseteq M(A\ot A).$$ 
In general, for a non-unital algebra, these two inclusions are strict. 
\snl
The opposite algebra $A^{\text{op}}$ is the algebra $A$, but endowed with the opposite product. 
\snl
We use $\iota$ for the identity map. We use $\zeta$ to denote the flip map on $A\ot A$ and for its extension to the multiplier algebra $M(A\ot A)$. A coproduct is a map from $A$ to $M(A\ot A)$. The composition of $\Delta$ with $\zeta$ is denoted by $\Delta^{\text{cop}}$. 
\snl
We will sometimes use the Sweedler notation for coproduct. So we write $\sum_{(a)} a_{(1)} \ot a_{(2)}$ for $\Delta(a)$ when $\Delta$ is a coproduct on the algebra $A$ and $a\in A$. Some useful information about the use of the Sweedler notation in the case of non-unital algebras can be found e.g.\ in [VD1].
\snl
For the theory of weak Hopf algebras we refer to [B-N-S] and for the relation with Hopf algebroids to [B].
\snl
For the theory of weak multiplier Hopf algebras, we refer to [VD-W1] and [VD-W2]. In particular, we find in [VD-W2] the necessary results about the source and target maps and source and target algebras that are very important for the described procedure. See also [VD-W0] for a motivational paper. The theory of multiplier Hopf algebroids has been developed recently in [T-VD]. An important role is played by separability idempotents. For these objects, we refer to [VD2].
\snl
The conventions used in the original papers on weak multiplier Hopf algebras [VD-W1] and [VD-W2] are not always the same as in the newer paper on multiplier Hopf algebroids [T-VD]. We will mention such a difference whenever it occurs. 
\nl

\bf Acknowledgements \rm
\nl
The second name author (Alfons Van Daele) likes to thank his coauthor for the hospitality during his visit to M\"unster where part of this work was done. He is also grateful to Piotr Hajac and the other organizers of the Focus Program on Noncommutative Geometry and Quantum Groups in Toronto (June 2013) for the invitation and the possibility to talk about this work at the Fields Institute during this program. Finally, we thank Gabriella B\"ohm for discussions on this subject.
\nl
\nl

\bf 1. Preliminaries \rm
\nl
In this preliminary section, we mainly recall the notion of a {\it regular weak multiplier Hopf algebra} as well as some of the basic properties. We also give the necessary definitions and results about the source and target maps and the source and target algebras.
\snl
We could have chosen to do the same here for the {\it multiplier Hopf algebroids}. However, it turned out to be more convenient to recall the notions and some of the results during the process in the next two sections, where we describe the passage from a regular weak multiplier Hopf algebra to a multiplier Hopf algebroid.
\snl
For details about the first objects, we refer to [VD-W1], [VD-W2] and [VD-W3]. The main reference for multiplier Hopf algebroids is [T-VD]. Further information can be found in the references given in the introduction.
\nl
\it Weak multiplier Hopf algebras \rm
\nl
The basic ingredient is a {\it non-degenerate} and {\it idempotent} algebra $A$ over the field of complex numbers and a {\it regular} and {\it full coproduct} $\Delta$ on $A$, admitting a {\it counit}.

\inspr{1.1} Definition \rm 
A {\it regular coproduct} on an algebra $A$ is a homomorphism \newline $\Delta:A\to M(A\ot A)$ so that elements of the form 
$$\align 
&\Delta(a)(1\ot b) \qquad\qquad\qquad (1\ot b)\Delta(a) \\
&\Delta(a)(b\ot 1) \qquad\qquad\qquad (b\ot 1)\Delta(a)
\endalign$$
all belong to $A\ot A$ for all $a,b\in A$.  It is assumed to satisfy {\it coassociativity}. This can be formulated in different ways, e.g. we assume
$$(a\ot 1\ot 1)(\Delta\ot\iota)(\Delta(b)(1\ot c))=(\iota\ot\Delta)((a\ot 1)\Delta(b))(1\ot 1\ot c)$$
for all $a,b,c\in A$. In any case, the first requirement is necessary to be able to formulate coassociativity in any of such forms. A coproduct is called {\it full} if the smallest subspaces $V,W$ of $A$ so that 
$$\Delta(a)(1\ot b)\in V\ot A
\qquad\quad\text{and}\qquad\quad
(a\ot 1)\Delta(b)\in A\ot W$$
for all $a,b\in A$ are actually all of $A$. A {\it counit} for the pair $(A,\Delta)$ is a linear map $\varepsilon:A\to \Bbb C$ satisfying
$$(\varepsilon\ot\iota)(\Delta(a)(1\ot b))=ab
\qquad\quad\text{and}\qquad\quad
(\iota\ot\varepsilon)((a\ot 1)\Delta(b))=ab$$
for all $a,b$.
\hfill$\square$
\einspr 

To formulate fullness of the coproduct, it is not needed that it is regular. The same is true for the characterizing properties of the counit. However we assume regularity of the coproduct from the very beginning. This is the only case that will be treated in this paper. For the more general notion, we refer to [VD-W1]. 
\snl
If the algebra has an identity, so that $\Delta:A\to A\ot A$, fullness of the coproduct follows from the existence of the counit. Also the uniqueness of the counit follows immediately from the definition in that case. However, as discussed in [VD-W0], this is not expected in general. On the other hand, for a full coproduct, it can be shown that the counit is unique. 
\snl
Also remark that, because the coproduct is assumed to be regular, the defining property of the counit can also be formulated with the factors on the other side. In fact, if $A$ is replaced by $A^{\text{op}}$, nothing above changes and $\Delta$ remains a regular and full coproduct and $\varepsilon$ is still the counit. Similarly if we replace $\Delta$ by $\Delta^{\text{cop}}$, obtained from $\Delta$ by composing it with the flip map $\zeta$ on $A\ot A$, extended to $M(A\ot A)$, nothing changes here.
\snl
A {\it regular weak multiplier Hopf algebra} is a pair $(A,\Delta)$ of a non-degenerate idempotent algebra with a regular and full coproduct, such that there is a counit and with some further assumptions, formulated in the terms of the ranges and the kernels of the {\it canonical maps} $T_1$, $T_2$, $T_3$ and $T_4$ given by 
$$\align 
T_1(a\ot b)&=\Delta(a)(1\ot b)
\qquad\quad\text{and}\quad\qquad
T_2(a\ot b)=(a\ot 1)\Delta(b)\\
T_3(a\ot b)&=(1\ot b)\Delta(a)
\qquad\quad\text{and}\quad\qquad
T_4(a\ot b)=\Delta(b)(a\ot 1)\endalign
$$
when $a,b\in A$. Remark that the maps $T_3$ and $T_4$ are the maps $T_1$ and $T_2$ for the opposite algebra $A^{\text{op}}$ (with the original coproduct).
\snl
For the {\it precise conditions} we refer to the original papers mentioned before.
\nl
We now collect the {\it main properties} that we will be using further in this paper.
\snl
In any weak multiplier Hopf algebra (regular or not), the following holds.

\inspr{1.2} Proposition \rm Let $(A,\Delta)$ be a weak multiplier Hopf algebra. There is a unique idempotent $E\in M(A\ot A)$ such that 
$$E(A\ot A)=\Delta(A)(A\ot A) 
\qquad\quad\text{and}\qquad\quad
(A\ot A)E=(A\ot A)\Delta(A).$$
It is the smallest idempotent element in $M(A\ot A)$ such that $\Delta(a)=E\Delta(a)=\Delta(a)E$ for all $a\in A$. The coproduct $\Delta$ can be extended uniquely to a homomorphism $\Delta:M(A)\to M(A\ot A)$ with the property that still $\Delta(m)=E\Delta(m)=\Delta(m)E$ for all $m\in M(A)$. Similarly we have extensions of the homomorphisms $\Delta\ot\iota$ and $\iota\ot\Delta$ to $M(A\ot A)$ and these are unique if we require that 
$$\align(\Delta\ot\iota)(m)&=(E\ot 1)((\Delta\ot\iota)(m))=((\Delta\ot\iota)(m))(E\ot 1)\\
(\iota\ot\Delta)(m)&=(1\ot E)((\iota\ot\Delta)(m))=((\iota\ot\Delta)(m))(1\ot E)
\endalign$$
for all $m\in M(A\ot A)$. For these extensions we have $\Delta(1)=E$ and also $(\Delta\ot\iota)E=(\iota\ot\Delta)E$. Finally also
$$(\Delta\ot\iota)E=(E\ot 1)(1\ot E).\tag"(1.1)"$$
\vskip -0.7 cm \hfill$\square$
\einspr

We are using the same symbols for the extensions of the homomorphisms involved. Also remark that the last property (1.1) is very natural. Indeed, by the construction of the extensions, we will have that $(\Delta\ot\iota)E$ is an idempotent, smaller than $E\ot 1$ and $1\ot E$. So the last property is saying that $(\Delta\ot\iota)E$ is as big as possible. Property (1.1) is sometimes referred to as the weak comultiplicativity of the unit. This makes sense once we know that $\Delta(1)=E$ for the extension of $\Delta$. See [K-VD] for the use of this terminology and [VD-W0] for a motivation of the equation (1.1).
\snl
One of the main results in the theory is the existence of a unique {\it antipode}. Recall that we assume {\it regularity} of the weak multiplier Hopf algebra. Then we have the following properties.

\inspr{1.3} Proposition \rm
There is a bijective linear map $S:A\to A$ such that the maps $R_1$ and $R_2$, defined on $A\ot A$ by
$$\align 
R_1(a\ot b)&=\sum_{(a)}a_{(1)}\ot S(a_{(2)})b \\
R_2(a\ot b)&=\sum_{(b)}aS(b_{(1)})\ot b_{(2)},
\endalign$$
are well-defined maps from $A\ot A$ to itself and so that $R_1$ and $R_2$ are generalized inverses of the canonical maps $T_1$ and $T_2$ respectively. 
\hfill $\square$
\einspr

We are using the Sweedler notation here. The equations can be covered by multiplying with an element of $A$ from the left in the first factor, for the first equation, and with an element of $A$ from the right in the second factor, for the second equation. It is not obvious that the right hand sides are in $A\ot A$ but that is actually the case. 
\snl
That $R_1$ is a generalized inverse of $T_1$ means that 
$$T_1R_1T_1=T_1 
\qquad\quad\text\quad\qquad
R_1T_1R_1=R_1$$
and similarly for $T_2$ and $R_2$.
\snl
In general, such generalized inverses are not unique. They are determined by a choice of the projection on the kernel and on the range of the corresponding maps. Such a choice is part of the axioms of a weak multiplier Hopf algebra and it depends only on the canonical multiplier $E$. It follows that also the map $S$, as in the above proposition, is uniquely determined. It is called the antipode of the weak multiplier Hopf algebra $(A,\Delta)$. As expected, it is an anti-isomorphism of the algebra $A$ and it also flips the coproduct $\Delta$.
\snl
If we replace $A$ by $A^{\text{op}}$, the antipode $S$ is replaced by the inverse $S^{-1}$. Then it is easy to find formulas for the appropriate generalized inverses $R_3$ and $R_4$ of $T_3$ and $T_4$ respectively.
\snl
From these equations, we can conclude the following result about the antipode.

\inspr{1.4} Proposition \rm
The antipode satisfies
$$\sum_{(a)} a_{(1)}S(a_{(2)})a_{(3)}=a
\qquad\quad\text{and}\quad\qquad
         \sum_{(a)} S(a_{(1)})a_{(2)}S(a_{(3)})=S(a)$$
for all $a\in A$. 
\hfill $\square$
\einspr

If we multiply the left hand side of the first equation with an element of $A$, left or right, it turns out that the use of the Sweedler notation is justified. Similarly for the second equation. So one first has to consider these equations in $M(A)$. As we get elements in $A$ on the right, we see that actually the equations will hold in $A$, but that is not trivial. We get the same formulas with the inverse $S^{-1}$ of the antipode, provided that we either flip the product or the coproduct. 
\snl
Because $R_1$ is a generalized inverse of $T_1$, it follows that $T_1R_1$ projects onto the range of $T_1$. Similarly for $T_2R_2$ and for the two other maps $T_3R_3$ and $T_4R_4$. It then follows from the choices mentioned earlier that the following holds.

\inspr{1.5} Proposition \rm
$$
\align T_1R_1(a\ot b)&=E(a\ot b)
\quad\qquad\text{and}\quad\qquad
T_2R_2(a\ot b)=(a\ot b)E\\
T_3R_3(a\ot b)&=(a\ot b)E
\quad\qquad\text{and}\quad\qquad
T_4R_4(a\ot b)=E(a\ot b)
\endalign$$
for all $a,b\in A$.
\hfill $\square$
\einspr

The second set of formulas is obtained from the first one by replacing $A$ by $A^{\text{op}}$ and taking into account that $E$ remains the same.
\snl
These formulas describe the ranges of these four canonical maps. This means that we have
$$
\align \Delta(A)(1\ot A)&=E(A\ot A)
\quad\qquad\text{and}\quad\qquad
(A\ot 1)\Delta(A)=(A\ot A)E\\
(1\ot A)\Delta(A)&=(A\ot A)E
\quad\qquad\text{and}\quad\qquad
\Delta(A)(A\ot 1)=E(A\ot A).
\endalign$$
In fact, these equations are part of the basic assumptions in the definition of a regular weak multiplier Hopf algebra.
\snl
We have similar formulas for the {\it kernels}.

\inspr{1.6} Proposition \rm
There exists idempotents  $F_1, F_2, F_3$ and $F_4$ in $M(A\ot A^{\text{op}})$ given by
$$\align F_1&=(\iota\ot S)E  \qquad\quad\text{and}\qquad\quad F_3=(\iota\ot S^{-1})E \\
	F_2&=(S\ot \iota)E \qquad\quad\text{and}\qquad\quad F_4=(S^{-1}\ot \iota)E.
\endalign$$
They give the kernels of the four canonical maps in the sense that 
$$\align R_1T_1(a\ot b)&=(a\ot 1)F_1(1\ot b)
\quad\quad\text{and}\quad\quad
R_2T_2(a\ot b)=(a\ot 1)F_2(1\ot b)\\
R_3T_3(a\ot b)&=(1\ot b)F_3(a\ot 1)
\quad\quad\text{and}\quad\quad
R_4T_4(a\ot b)=(1\ot b)F_4(a\ot 1)
\endalign$$
for all $a,b$ in $A$.
\hfill $\square$ 
\einspr

These results follow from the formulas in the previous proposition, provided we use that $S$ is an anti-isomorphism of the algebra $A$. Indeed, we have e.g.
$$R_1(\iota\ot S)=(\iota\ot S)T_3
\qquad\quad\text{and}\qquad\quad
T_1(\iota\ot S)=(\iota\ot S)R_3$$
and this will give the formula for $R_1T_1$. Similarly for the other cases.
\snl
From these formulas, we also find
$$\align \text{Ker}(T_1)&=(A\ot 1)(1-F_1)(1\ot A)
\quad\quad\text{and}\quad\quad
\text{Ker}(T_2)=(A\ot 1)(1-F_2)(1\ot A)\\
\text{Ker}(T_3)&=(1\ot A)(1-F_3)(A\ot 1)
\quad\quad\text{and}\quad\quad
\text{Ker}(T_4)=(1\ot A)(1-F_4)(A\ot 1).
\endalign$$
We refer to [VD-W1] for these results. In that paper, it is also shown that the multipliers $F_i$ are completed determined by $E$ itself and the antipode is not needed to characterize them.
\nl
\it The source and target maps \rm
\nl
Having formulated the basic results about regular weak multiplier Hopf algebras, we now recall the {\it source} and {\it target maps}, the {\it source and target algebras} and their main properties needed further in this paper.

\inspr{1.7} Definition \rm The source and target maps $\varepsilon_s$ and $\varepsilon_t$ from $A$ to $M(A)$ are defined by
$$\varepsilon_s(a)=\sum_{(a)}S(a_{(1)})a_{(2)}
\qquad\quad\text{and}\quad\qquad
\varepsilon_t(a)=\sum_{(a)}a_{(1)}S(a_{(2)}).
$$
\vskip -0.8 cm
\hfill $\square$
\einspr

Recall that by assumption, $\varepsilon_s(a)$ as defined above, is a left multiplier. It is also a right multiplier because of the result in Proposition 1.3. Similarly $\varepsilon_t(a)$ is a right multiplier by the assumptions on the coproduct and a left multiplier because of Proposition 1.3.
\snl
In what follows, we denote $\varepsilon_s(A)$ by $B$ and $\varepsilon_t(A)$ by $C$. 
\snl
For any $a\in A$, $x\in B$ and $y\in C$ we have
$$\align
\varepsilon_s(ax)&=\varepsilon_s(a)x
\qquad\quad\text{and}\qquad\quad
\varepsilon_s(ay)=S(y)\varepsilon_s(a)\tag"(1.2)"\\
\varepsilon_t(ya)&=y\varepsilon_t(a)
\qquad\quad\text{and}\qquad\quad
\varepsilon_t(xa)=\varepsilon_t(a)S(x).\tag"(1.3)"
\endalign$$

\inspr{1.8} Proposition \rm
The images $B$ and $C$ of the source and target maps $\varepsilon_s(A)$ and $\varepsilon_t(A)$ are non-degenerate idempotent algebras.  They sit nicely in the multiplier algebra $M(A)$ in the sense that 
$$BA=AB=A
\qquad\quad\text{and}\qquad\quad
CA=AC=A.$$
As a consequence, the imbeddings extend to imbeddings of their multiplier algebras $M(B)$ and $M(C)$ in $M(A)$. These multiplier algebras are denoted by $A_s$ and $A_t$ respectively. They are called the {\it source} and {\it target algebras} and they are characterized in $M(A)$ by
$$\align A_s&=\{x\in M(A)\mid \Delta(x)=E(1\ot x)\}\\
A_t&=\{y\in M(A)\mid \Delta(y)=(y\ot 1)E\}.
\endalign$$
\vskip -0.8cm
\hfill $\square$
\einspr

The formulas in (1.2) and (1.3) also hold for $x\in A_s$ and $y\in A_t$.
\snl
The ranges of the source and target maps are, in a sense that can be made precise, the left and the right leg of $E$ respectively. Consequently the algebras $A_s$ and $A_t$ commute. This implies that the formulas in the proposition above still hold when the product in $A$ is reversed.
\snl
Finally we have the following important property.

\inspr{1.9} Proposition \rm
We have $E\in M(B\ot C)$ and $E$ is a separability idempotent in the sense of [VD2]. The associated antipodal maps from $B$ to $C$ and from $C$ to $B$ are precisely the restriction of the antipode (or rather of its extension to the multiplier algebra $M(A)$). This means that
$$E(x\ot 1)=E(1\ot S(x))
\qquad\quad\text{and}\quad\qquad
(1\ot y)E=(S(y)\ot 1)E
$$
for all $x\in B$ and $y\in C$.
\hfill $\square$
\einspr

Remark that $E(1\ot y)\in B\ot C$ for all $y\in C$ and if we use a Sweedler type notation $E_1\ot E_2$ for $E$ we have $S_B(E_1)E_2y=y$ when $y\in C$. Similarly, $(x\ot 1)E\in B\ot C$ and $xE_1S_C(E_2)=x$ for all $x\in B$. We write these properties as
$$S_B(E_1)E_2=1
\qquad\quad\text{and}\qquad\quad
E_1S_C(E_2)=1.$$
See Proposition 1.8 in [VD2].
\snl
We finish this section with the remark that the algebras $B$ and $C$ have local units (see Proposition 1.9 in [VD2]). The same is true for the original algebra $A$ (see Proposition 4.9 in [VD-W1]).
\nl
\nl

\bf 2. Compatible quantum graphs and separability idempotents \rm
\nl
In this section, we investigate the relation between compatible quantum graphs and separability idempotents. {\it Compatible quantum graphs}, as studied here, have been introduced in [T-VD] and are the basic ingredients for a multiplier Hopf algebroid. On the other side, {\it separability idempotents}, as studied in [VD2], play a basic role in the theory of weak multiplier Hopf algebras. And since this paper is devoted to the relation between these two objects, it is quite natural to start with the study of the relation between these basic ingredients.
\snl
Observe that in this section, {\it we do not yet consider the coproducts}. This will be done in the next two sections.
\nl
In a first item of this section, we will collect properties about separability idempotents (as introduced in [VD2]) and the concept of a separable Frobenius algebra.

\newpage
\it Separable Frobenius algebras \rm
\nl
We consider an algebra $B$ with local units. In particular, the product is non-degenerate and $B$ is idempotent. 
\snl
A linear functional $\varphi$ on $B$ is called {\it faithful} if the maps
$$x\mapsto \varphi(\,\cdot\,x)
\qquad\quad\text{and}\qquad\quad
x\mapsto\varphi(x\,\cdot\,)$$
from $B$ to the linear dual $B'$ are both injective. If these maps also have the same range, there exists a unique automorphism $\sigma$ of $B$ so that 
$\varphi(x_1x_2)=\varphi(x_2\sigma(x_1))$ for all $x_1$, $x_2$ in $B$. Because $B$ is assumed to be idempotent, this implies also that $\varphi$ is $\sigma$-invariant.
\snl
We call $\sigma$ the {\it modular automorphism} of $\varphi$.
\snl
Observe that the existence of a faithful linear functional on the algebra $B$ automatically implies that the product is non-degenerate.
\snl
If $B$ is a finite-dimensional algebra with identity, the existence of a faithful linear functional is equivalent with $B$ being Frobenius. Moreover, in this case, a faithful linear functional will automatically have a modular automorphism. In the literature, the inverse of $\sigma$ is called the Nakayama automorphism. We refer to [M] for the notion of the Nakayama automorphism in the case of (finite-dimensional) Frobenius algebras. However, in the case of an infinite-dimensional algebra, possibly without identity, this result is not true. There can be a linear functional that is faithful but without a modular automorphism. We refer to [VD2] for an example.
\nl
We now recall the following concept from [VD2].

\inspr{2.1} Definition \rm 
Let $\varphi$ be a linear functional on the algebra $B$. Assume that it is faithful and admits a modular automorphism as above. Let $C$ be the algebra $B^{\text{op}}$ and use $S$ for the identity map from $B$ to $C$. It is an anti-isomorphism by definition. Assume that there is an idempotent $E\in M(B\ot C)$ satisfying 
$$(\varphi(\,\cdot\,x)\ot \iota)E=S(x)\tag"(2.1)"$$
for all $x$. Then we call $\varphi$ a {\it separating functional} on $B$. If the algebra $B$ admits such a separating functional, we call $B$ {\it separable Frobenius}.
\hfill$\square$
\einspr

We need to give a couple of remarks here.

\inspr{2.2} Remarks \rm
i) If $B$ is a finite-dimensional algebra with identity and $\varphi$ a faithful linear functional, there is always an element $E\in B\ot C$ determined by the equality 2.1. However, it can happen that $E$ is not an idempotent. In fact one can easily construct examples where $E^2=0$. See [VD2] for an example.
\snl
ii) Given $\varphi$ the idempotent $E$ is unique if it exist. This follows from the faithfulness of $\varphi$. This means that we really define a property of the linear functional $\varphi$ on the algebra $B$.
\snl
iii) It can be shown that, if a separating functional on the algebra exists, then any faithful linear functional will have a modular automorphism and it will be separating as well. 
\hfill$\square$
\einspr

In what follows, we assume that $B$ is separable Frobenius and that $\varphi$ is a given separating functional on $B$ as in Definition 2.1. The modular automorphism (that exists by assumption) will be denoted by $\sigma_B$. Finally, the given anti-isomorphism from $B$ to $C$ will now be denoted by $S_B$.
\snl
Then we have the following result. 

\inspr{2.3} Proposition \rm We have an anti-isomorphism $S_C:C\to B$ defined by $S_C=\sigma_B^{-1}S_B^{-1}$. Define a linear functional $\varphi_C$ on $C$ by $\varphi_C(y)=\varphi_B(S_B^{-1}(y))$ when $y\in C$. It is a faithful functional with modular automorphism $\sigma_C$ given by $\sigma_C=S_BS_C$. We not only have $(\varphi(\,\cdot\,x)\ot \iota)E=S_B(x)$ for all $x\in B$ but now also
$$(\iota\ot\varphi_C(y\,\cdot\,))E=S_C(y)$$
for all $y\in C$. The element $E$ is a separability idempotent in $M(B\ot C)$ with associated anti-isomorphisms $S_B$ and $S_C$. This means that 
$$E(x\ot 1)=E(1\ot S_B(x))
\qquad\quad\text{and}\qquad\quad
(1\ot y)E=(S_C(y)\ot 1)E$$
for all $x\in B$ and $y\in C$. The linear functionals $\varphi_B$ and $\varphi_C$ are the integrals.
\hfill $\square$ \einspr

The proof is straightforward using techniques from [VD2].
\snl
So we see that given a separable Frobenius algebra $B$ with a separating linear functional $\varphi$, the unique associated idempotent $E$ is a separability idempotent in $M(B\ot C)$ where $C=B^{\text{op}}$. The identity map from $B$ to $C$ is the antipodal map $S_B:B\to C$, characterized by the equation $E(x\ot 1)=E(1\ot S_B(x))$ for all $x\in B$. Remark that it is not essential that $C$ is identified with $B^{\text{op}}$. Obviously, any algebra $C$ with a given anti-isomorphism $S_B:B\to C$ as a starting point will work as well.
\snl
In fact, also if we have two faithful functionals $\varphi_B$ and $\varphi'_B$, related by an automorphism of $B$, we will still get essentially the same separability idempotent. 
\snl
We also have the following converse result.

\inspr{2.4} Proposition \rm Suppose that $B$ and $C$ are non-degenerate algebras and that $E$ is a separability idempotent in $M(B\ot C)$. Then $B$ is separable Frobenius. The unique linear functional $\varphi_B$ on $B$ satisfying $(\varphi_B\ot\iota)E=1$ in $M(C)$ is a separating functional on $B$. 
\hfill $\square$ \einspr

In Section 4, where we characterize those multiplier Hopf algebroids that arise from a weak multiplier Hopf algebra, we will start with a multiplier Hopf algebroid whose base algebra $B$ is separable Frobenius as defined in Definition 2.1 above. The separability idempotent $E$ that we find in Proposition 2.3 will eventually be the canonical idempotent of the associated weak multiplier Hopf algebra. In Section 3, we start from a weak multiplier Hopf algebra and construct a multiplier Hopf algebroid from it. In this case we have given the canonical idempotent $E$. It is a separability idempotent. 
\nl
Next we consider the {\it second basic ingredient}, used in the theory of multiplier Hopf algebroids. First we look at it independently from the previous item. Later in this section, we will relate the two concepts.
\nl
\it Compatible quantum graphs \rm
\nl
We start with an algebra $A$ and again we assume that it has local units. Then it is automatically non-degenerate and idempotent. We also assume that $B$ and $C$ are two commuting subalgebras of $M(A)$ satisfying  
$$BA=AB=A
\qquad\quad\text{and}\quad\qquad
CA=AC=A.$$
It follows from these conditions that $B$ and $C$ are also non-degenerate algebras and that the imbeddings of $B$ and $C$ in $M(A)$ extend in a unique way to imbeddings of the multiplier algebras $M(B)$ and $M(C)$ in $M(A)$. In fact we have the following characterization
$$\align
M(B)&=\{m\in M(A) \mid mx, xm\in B \text{ for all } x\in B\} \\
M(C)&=\{m\in M(A) \mid my, ym\in C \text{ for all } y\in C\}.
\endalign$$
See e.g.\ Lemma 2.11 in [VD-W2] for an argument.
\snl
We denote elements of $B$ with $x,x',...$ and elements of $C$ with $y,y',...$ as is done in [T-VD]. We reserve $a,b,c$ for elements of $A$.
\snl
We also assume the existence of two anti-isomorphisms
$$ S_B:B\to C 
\qquad\quad\text{and}\quad\qquad
S_C:C\to B.
$$
As mentioned before, these two anti-isomorphisms are in general not expected to be each others inverses.

\inspr{2.5} Proposition \rm With the above assumptions, the two quadruples
$$
\Cal A_B:=(B,A,\iota_B, S_B)
\qquad\quad\text{and}\quad\qquad
\Cal A_C:=(C,A,\iota_C,S_C)
$$
are a compatible pair of a left and right quantum graph.

\snl\bf Proof\rm:
i) First we need to argue that $\Cal A_B$, composed of the ingredients $(B,A,\iota_B,S_B)$, is a left quantum graph as introduced in Definition 2.1 of [T-VD]. All of this is simple and quite straightforward. Remark that we have assumed from the very beginning that $A$ has local units. Then condition (4) of Definition 2.1 of [T-VD] is fulfilled.
\snl
ii) In a completely similar way $\Cal A_C$, composed of $(C,A,\iota_C,S_C)$, is a right quantum graph as in Definition 5.1 of [T-VD].
\snl
iii) Finally one can easily check that the two are compatible in the sense of Definition 6.1 of [T-VD].
\hfill$\square$
\einspr

It is more or less obvious from what we have seen in the preliminary section, that a regular weak multiplier Hopf algebra $(A,\Delta)$ gives rise in a natural way to such a compatible pair of quantum graphs. In this case, $B$ and $C$ are the images of the source and target maps, as sitting in $M(A)$, while $S_B$ and $S_C$ are the restrictions of the antipode (after is has been extended to the multiplier algebra). This will be the starting point in the next section where we pass from a regular weak multiplier Hopf algebra to a regular multiplier Hopf algebroid.
\snl
In Section 4, we will start from the data as in this section.
\snl
Passing again to the general situation, we consider several balanced tensor products. We can e.g.\ consider the quotient of $A\ot A$ over the subspace spanned by elements of the form $xa\ot b-a\ot S_B(x)b$ with $a,b\in A$ and $x\in B$. This space will be denoted as $A\ot _\ell A$. 
\snl
There are all together six cases that we will consider. We can define them in the present setting, but since we will only work further with these spaces in the third item of this section, we postpone the definitions of the five other cases. See Notation 2.7 below.
\snl
Indeed, next we combine the results from to two preceding ones.
\nl
\it Compatible quantum graphs with a separable Frobenius base algebra \rm
\nl
Again start with an algebra $A$ with local units, two commuting subalgebras $B$ and $C$ sitting nicely in $M(A)$ and two anti-isomorphisms $S_B:B\to C$ and $S_C:C\to B$ as in the previous item of this section. In other words, we assume that the quadruples 
$$
\Cal A_B:=(B,A,\iota_B, S_B)
\qquad\quad\text{and}\quad\qquad
\Cal A_C:=(C,A,\iota_C,S_C)
$$
are a compatible pair of a left and right quantum graph as defined in [T-VD]. 
\snl
Now we {\it add the following assumption} thus making a link between the material studied in the first item and the one discussed in the second item of this section.

\inspr{2.6} Assumption \rm We assume that $B$ (and hence also $C$) is separable Frobenius as in Definition 2.1. In addition, we require the existence of a separating linear functional $\varphi_B$ as in Definition 2.1 with the {\it extra assumption} that the modular automorphism $\sigma_B$ of $\varphi_B$ coincides with the inverse of the given automorphism $S_CS_B$ of $B$.
\hfill$\square$\einspr

It follows that we have a separability idempotent $E\in M(B\ot C)$ so that $S_B$ and $S_C$ are the antipodal maps associated with $E$ and that $\varphi_B$ is the left integral as in Proposition 2.1 of [VD2]. 
\snl
In other words, we have two commuting subalgebras $B$ and $C$ sitting nicely in $M(A)$ and a separability idempotent $E\in M(B\ot C)$. This is an {\it equivalent way} to formulate the setting in this item of the present section.
\snl
In Section 3 we arrive at this situation when we start from a regular weak multiplier Hopf algebra. In Section 4, it will be part of the assumptions.
\nl
Now we will consider the {\it various balanced tensor products} as announced already at the end of the second item of this section. Observe again that these spaces can be considered for any compatible pair of quantum graphs. It is not needed that the base algebra is separable Frobenius for this.

\inspr{2.7} Notation \rm We use the following notations:
\snl
\parindent 1.5cm
\item{(1.1)} \quad $A\ot _\ell A$ \ for the quotient of $A\ot A$ defined by $xa\ot b=a\ot S_B(x)b$ and
\item{(1.2)} \quad $A\ot _r A$ \ for the quotient of $A\ot A$ defined by $a\ot by=S_C(y)a\ot b$
\inspr{}\rm
where $x\in B$ and $y\in C$ and of course $a,b\in A$. Similarly we use
\parindent 1.5 cm \snl
\item{(2.1)} \quad $A\ot_s A$ \ defined by $ax\ot b=a\ot xb$,
\item{(2.2)} \quad $A\ot_t A$ \ defined by $a\ot yb=ay\ot b$,
\item{(2.3)} \quad $A\ot^s A$ \ defined by $xa\ot b=a\ot bx$ and
\item{(2.4)} \quad $A\ot^t A$ \ defined by $a\ot by=ya\ot b$,
\inspr{} \rm
where again $x\in B$, $y\in C$ and $a,b\in A$. 
\hfill$\square$\einspr

We will now obtain {\it concrete realizations} of these balanced tensor products in terms of the separability idempotent $E$. We begin with the first two cases, namely $A\ot_\ell A$ and $A\ot_r A$.

\inspr{2.8} Proposition \rm
i) Denote by $\pi_\ell$ the canonical quotient map from $A\ot A$ to $A\ot_\ell A$. There is a well-defined map $\theta_\ell:A\ot_\ell A\to A\ot A$ given by $\theta_\ell(a\ot b)=E(a\ot b)$. It satisfies $\pi_\ell\circ \theta_\ell=\iota_\ell$, where $\iota_\ell$ denotes the identity map on $A\ot_\ell A$. 
\snl
ii) Denote by $\pi_r$ the canonical quotient map from $A\ot A$ to $A\ot_r A$. There is a well-defined map $\theta_r:A\ot_r A\to A\ot A$ given by $\theta_r(a\ot b)=(a\ot b)E$. It satisfies $\pi_r\circ \theta_r=\iota_r$, where $\iota_r$ denotes the identity map on $A\ot_r A$. 

\snl\bf Proof\rm:
i) For all elements $a,b\in A$ and $x\in B$ we have 
$$E(xa\ot b)=E(a\ot S_B(x)b)$$
and therefore the map $a\ot b\mapsto E(a\ot b)$ is a well-defined map $\theta_\ell$ from the balanced tensor product $A\ot_\ell A$ to $A\ot A$. Furthermore, for all elements $a,b\in A$ and $y\in C$, we have
$$
\align \pi_\ell(\theta_\ell(a\ot yb)&=\pi_\ell(E(a\ot yb))\\
&=\pi_\ell(a\ot S_B(E_{(1)})E_{(2)}yb)\\
&=\pi_\ell(a\ot yb)
\endalign$$
where we are using the Sweedler type notation $E=E_{(1)}\ot E_{(2)}$. Remark that 
$$E(1\ot C)\subseteq B\ot C
\qquad\quad\text{and}\qquad\quad
S_B(E_{(1)})E_{(2)}y=y$$ 
when $y\in C$ as shown in Proposition 1.8 of [VD2] (see also a remark at the end of the preliminary section). 
Finally, because $CA=A$, we get $\pi_\ell(\theta_\ell(a\ot b)=a\ot b$ in $A\ot_\ell A$ for all $a,b\in A$. 
\snl
ii) The proof here is completely similar. Now we use that 
$$(1\ot y)E=(S_C(y)\ot 1)E$$ 
for all $y\in C$ and that $xE_{(1)}S_C(E_{(2)})=x$ for all $x\in B$ (again see Proposition 1.8 of [VD2] and a remark at the end of the preliminary section). 
\hfill$\square$\einspr
 
It follows from this result that the map $\theta_\ell:A\ot_\ell A\to A\ot A$ is injective and its range is $E(A\ot A)$. Similarly the map $\theta_r:A\ot_r A\to A\ot A$ is injective with range $(A\ot A)E$.
\snl
Now we prove similar results for the four other balanced products in Notation 2.8. First we need the following analogue of (part of) the result reviewed in Proposition 1.6 in the case of a weak multiplier Hopf algebra.

\inspr{2.9} Proposition \rm
There exist idempotents  $F_1, F_2, F_3$ and $F_4$ in $M(A\ot A^{\text{op}})$ given by
$$\align F_1&=(\iota\ot S_C)E  \qquad\quad\text{and}\qquad\quad F_3=(\iota\ot S_B^{-1})E \\
	F_2&=(S_B\ot \iota)E \qquad\quad\text{and}\qquad\quad F_4=(S_C^{-1}\ot \iota)E.
\endalign$$

\snl\bf Proof\rm:
Recall that $E(x\ot 1)\in B\ot C$ for all $x\in B$. Then $E(xa\ot 1)\in A\ot C$ and as $BA=A$ we find $E(A\ot 1)\subseteq A\ot C$. It follows that $(\iota\ot S_C)E$ is a well-defined multiplier $F_1$ in $M(A\ot A^{\text{op}})$. Similarly for the other elements $F_2$, $F_3$ and $F_4$. 
\hfill $\square$ 
\einspr 

Remark that the identity map  from $A\ot A$ to itself is an anti-isomorphism from $A\ot A^{\text{op}}$ to $A^{\text{op}}\ot A$ that extends to a natural anti-isomorphism from $M(A\ot A^{\text{op}})$ to $M(A^{\text{op}}\ot A)$. 
\snl
Now we find the following results about the four other balanced tensor products.

\iinspr{2.10} Proposition \rm
We can define the four maps
$$\align
&a\ot b\mapsto (a\ot 1)F_1(1\ot b) \qquad \text{ on } \qquad A\ot_s A \\
&a\ot b\mapsto (a\ot 1)F_2(1\ot b) \qquad \text{ on } \qquad A\ot_t A \\
&a\ot b\mapsto (1\ot b)F_3(a\ot 1) \qquad \text{ on } \qquad A\ot^s A \\
&a\ot b\mapsto (1\ot b)F_4(a\ot 1) \qquad \text{ on } \qquad A\ot^t A, 
\endalign$$
all of them with range in $A\ot A$. For each of them, when composed with the corresponding canonical projection map, we find the identity map on the corresponding balanced tensor product.

\snl\bf Proof\rm:
In the first case, we  use e.g.\ that 
$$(x\ot 1)F_1=(1\ot S_C)((x\ot 1)E)\in B\ot B$$
and that $xE_{(1)}S_C(E_{(2)})=x$ for all $x\in B$. In the second case we use that 
$$F_2(1\ot y)=(S_B\ot \iota)(E(1\ot y))\in C\ot C$$ 
and that $S_B(E_{(1)})E_{(2)}y=y$ for all $y\in C$. For the third case, we need 
$$F_3(x\ot 1)=(\iota\ot S_B^{-1})(E(x\ot 1))\in B\ot B$$ 
and $S_B^{-1}(E_{(2)})E_{(1)}x=x$ for all $x\in B$. Finally in the last case we use
$$(1\ot y)F_4=(S_C^{-1}\ot \iota)((1\ot y)E)\in C\ot C$$ and 
$yE_{(2)}S_C^{-1}(E_{(1)})=y$ for all $y\in C$. 
\hfill$\square$\einspr

Again, as a consequence these four maps are bijections of these balanced tensor product with appropriate subspaces of $A\ot A$.
\nl\nl

\bf 3. From weak multiplier Hopf algebras to multiplier Hopf algebroids \rm
\nl
In this section we start with a regular weak multiplier Hopf algebra $(A,\Delta)$. We set 
$$B=\varepsilon_s(A)
\qquad\quad\text{and}\qquad\quad
C=\varepsilon_t(A)$$
where $\varepsilon_s$ and $\varepsilon_t$ are the source and target maps. We know that $B$ and $C$ are non-degenerate idempotent algebras, sitting nicely in the multiplier algebra $M(A)$ of $A$. They are commuting subalgebras of $M(A)$. The canonical idempotent $E$ of the weak multiplier Hopf algebra is a separability idempotent in $M(B\ot C)$. The associated antipodal maps $S_B:B\to C$ and $S_C:C\to B$ are the corresponding restrictions of the extension of the antipode $S$ to the multiplier algebra $M(A)$. Because of this, it can not be confusing if we sometimes drop the indices here and simply use $S$ for these two associated maps. We refer to the preliminary section for more details.
\snl 
We have seen in the previous section that the regular weak multiplier Hopf algebra $(A,\Delta)$ gives rise in a canonical way to a pair of compatible quantum graphs
$$
\Cal A_B:=(B,A,\iota_B, S_B)
\qquad\quad\text{and}\quad\qquad
\Cal A_C:=(C,A,\iota_C,S_C)
$$
as shown in Proposition 2.5.
\snl
The {\it aim of this section} is to show that the coproduct $\Delta$ induces a left and a right coproduct on $A$ making it into a regular multiplier Hopf algebroid as defined in [T-VD]. We will also express the data of the multiplier Hopf algebroid in terms of the various data of the given weak multiplier Hopf algebra. 
\snl
In the {\it next section}, we find a necessary and sufficient condition for a multiplier Hopf algebroid to arise in this way.
\nl
\it The left and the right coproducts $\Delta_B$ and $\Delta_C$ and the main result \rm
\nl
We first consider the balanced tensor products $A\ot_\ell A$ and $A\ot_r A$ as in Notation 2.7 of the previous section.
\snl 
The first space $A\ot_\ell A$ is endowed with the two actions of $A$, obtained by multiplication from the right in the first and in the second factor. We use $A\overline\ot_\ell A$ for the {\it extended module}. It can be characterized by the property that for any $z\in A\overline\ot_\ell A$, we have well-defined elements 
$$z(1\ot a)
\quad\qquad\text{and}\quad\qquad
z(a\ot 1)$$
in $A\ot_\ell A$ for all $a\in A$. Similarly we use $A\overline\ot_r A$ for the extended module of $A\ot_r A$ where elements $z\in A\overline\ot_r A$ are characterized by the property that 
$$(a\ot 1)z
\quad\qquad\text{and}\quad\qquad
(1\ot a)z$$
are in $A\ot_r A$ for all $a\in A$. 
\snl
For a more detailed treatment of the notion of an extended module, sometimes also called the completed module, we refer to Section 2 of [VD1].

\inspr{3.1} Notation \rm
We use $L_{\text{reg}}(A_B\times A)$ for the subspace of $A\overline\ot_\ell A$ of elements $z$  with the property that 
$$z(x\ot 1)=z(1\ot S(x))$$
for all $x\in B$. Remark that we extend the module action to the multiplier algebras. Similarly, we use $R_{\text{reg}}(A\times_C A)$ for the subspace of elements $z\in A\overline\ot_r A$ with the property that 
$$(1\ot y)z=(S(y)\ot 1)z$$
for all $y\in C$. Again we extend the module action to the multiplier algebras. 
\hfill $\square$ \einspr

There is a natural imbedding of $L_{\text{reg}}(A_B\times A)$ into the space $\text{End}(A\ot_\ell A)$ by letting an element $\xi$ of $L_{\text{reg}}(A_B\times A)$ act on $A\ot_\ell A$ simply by 
$$\xi\pi_\ell(a\ot b)=\xi(a\ot b)$$
for $a,b\in A$ where as before, $\pi_\ell$ denotes the canonical projection of $A\ot A$ onto $A\ot_\ell A$. In the right hand side, $\xi(a\ot b)$ is the result of the right action of $a\ot b$ on $\xi$. The action of $L_{\text{reg}}(A_B\times A)$ on $A\ot_\ell A$ is well-defined, precisely because of the condition imposed on elements in $L_{\text{reg}}(A_B\times A)$. It is also not hard to see that composition of maps induces a product on $L_{\text{reg}}(A_B\times A)$, making it into an associative algebra. 
\snl
The algebra $L_{\text{reg}}(A_B\times A)$ is defined in Definition 2.4 of [T-VD] and is called the algebra of regular left multipliers of the left Takeuchi product $A_B\times A$. Similarly the space \newline $R_{\text{reg}}(A\times_C A)$ is defined in Definition 5.3 of [T-VD] where it is called the algebra of regular right multipliers of the right Takeuchi product $A\times_C A$. 
\snl
These two spaces are the target spaces for the left and the right coproducts respectively in what follows.

\inspr{3.2} Proposition \rm Define $\Delta_B:A\to  A\overline\ot_\ell A$ by
$$\Delta_B(a)(1\ot b)=\pi_\ell(\Delta(a)(1\ot b))
\quad\quad\text{and}\quad\quad
\Delta_B(a)(c\ot 1)=\pi_\ell(\Delta(a)(c\ot 1))$$ 
for $a,b,c\in A$. Here again $\pi_\ell$ is the canonical projection of $A\ot A$ onto the quotient space $A\ot_\ell A$. Then $\Delta_B$ has range in $L_{\text{reg}}(A_B\times A)$ and it is a left coproduct on the left quantum graph $(B,A,\iota, S_B)$ in the sense of Definition 2.11 of [T-VD]. 
\snl
Similarly, define $\Delta_C:A\to A\overline\ot_r A$ by 
$$(c\ot 1)\Delta_C(a)=\pi_r((c\ot 1)\Delta(a))
\quad\quad\text{and}\quad\quad
(1\ot b)\Delta_C(a)=\pi_r((1\ot b)\Delta(a))$$
where now we use $\pi_r$ for the canonical projection of $A\ot A$ onto the quotient space $A\ot_r A$. Then $\Delta_C$ has range in $R_{\text{reg}}(A\times_C A)$ and it is a right coproduct on the right quantum graph $(C,A,\iota, S_C)$ in the sense of Definition 5.6 of [T-VD]. 
\snl
The two coproducts are compatible in the sense of Definition 6.4 of [T-VD].

\snl\bf Proof\rm:
The proof is rather straightforward. The different steps are as follows.
\snl
i) First observe that $\Delta_B(a)$ is a well-defined element in the extended module $A\overline\ot_\ell A$. To show this, first the elements 
$\Delta_B(a)(1\ot b)$ and $\Delta_B(a)(c\ot 1)$, given by the formulas in the formulation, are considered and then the obvious compatibility relations are verified. All of this is essentially trivial.
\snl
ii) Next it has to be shown that $\Delta_B(a)$, defined in this way, really belongs to the subspace $L_{\text{reg}}(A_B\times A)$. This follows from the fact that
$$\Delta(a)(x\ot 1)=\Delta(a)(1\ot S(x))$$
for all $a\in A$ and $x\in B$. 
\snl
iii) The map $\Delta_B$ will be a homomorphism from $A$ to the algebra $L_{\text{reg}}(A_B\times A)$ because $\Delta$ is a homomorphism. 
\snl
iv) The behavior of $\Delta_B$ on $B$ and $C$ is determined by the behavior of $\Delta$ itself on these subalgebras.
\snl
v) Finally, coassociativity of $\Delta_B$ will follow from coassociativity of $\Delta$, written in the form
$$(\Delta\ot\iota)(\Delta(a)(1\ot b))(c\ot 1\ot 1)=(\iota\ot\Delta)(\Delta(a)(c\ot 1))(1\ot 1\ot b)$$
for all $a,b,c\in A$. We have to project onto the appropriate balanced subspace of $A\ot A\ot A$ and use that 
$$\Delta(xa)=(1\ot x)\Delta(a)
\qquad\quad\text{and}\qquad\quad
\Delta(ya)=(y\ot 1)\Delta(a)$$
when $a\in A$, $x\in B$ and $y\in C$.
\snl
The same arguments are used to show that $\Delta_C$ is a right coproduct.
\snl
The compatibility of the two coproducts is again obtained from the coassociativity of $\Delta$, now formulated as
$$\align
(c\ot 1\ot 1)(\Delta\ot\iota)(\Delta(a)(1\ot b))
&=(\iota\ot\Delta)(c\ot 1)\Delta(a))(1\ot 1\ot b)\\
(1\ot 1\ot b)(\iota\ot\Delta)(\Delta(a)(c\ot 1))
&=(\Delta\ot\iota)((1\ot b)\Delta(a))(c\ot 1\ot 1)
\endalign$$
for all $a,b,c\in A$. Also here, we have to project onto the appropriate quotient space of $A\ot A\ot A$.
\hfill $\square$
\einspr

Next we consider the {\it canonical linear maps} associated with these coproducts. Let us first recall the definitions.

\inspr{3.3} Notation \rm For the left coproduct $\Delta_B$ we have the maps $T_\lambda$ and $T_\rho$ going from $A\ot A$ to $A\ot_\ell A$ and defined by
$$T_\lambda(a\ot b)=\Delta_B(b)(a\ot 1)
\qquad\quad\text\qquad\quad
T_\rho(a\ot b)=\Delta_B(a)(1\ot b)$$
for $a,b,c\in A$; see Section 3 of [T-VD]. For the right coproduct $\Delta_C$ we have the maps $_\lambda T$ and $_\rho T$, going from $A\ot A$ to $A\ot_r A$ and defined by
$$_\lambda T(a\ot b)=(a\ot 1)\Delta_C(b)
 \qquad\quad\text\qquad\quad
 _\rho T(a\ot b)=(1\ot b)\Delta_C(a)$$
for $a,b,c\in A$; see Section 5 of [T-VD].
\einspr
 
It is useful to compare these four maps with the canonical maps $T_1$, $T_2$, $T_3$ and $T_4$ as used in [VD-W1]. These are maps from $A\ot A$ to itself, given by the formulas 
$$\align T_1(a\ot b)&=\Delta(a)(1\ot b)
	\qquad\quad\text{and}\qquad\quad
		T_2(a\ot b)=(a\ot 1)\Delta(b)\\
T_3(a\ot b)&=(1\ot b)\Delta(a)
	\qquad\quad\text{and}\qquad\quad
		T_4(a\ot b)=\Delta(b)(a\ot 1).
\endalign$$
See Section 1 in [VD-W1]. We get from one set to the other by applying the appropriate quotient maps, sometimes combined with the flip. 
\snl
As a consequence of the results in the previous section on the various balanced tensor products (Proposition 2.8 and Proposition 2.10), we find the following. 

\inspr{3.4} Proposition \rm
The maps, as denoted in Notation 3.3, give rise to {\it bijective maps}, still denoted (for simplicity) with the same symbols,
$$\align &T_\lambda:A\ot^t A\to A\ot_\ell A
\qquad\quad\text{and}\qquad\quad
T_\rho:A\ot_s A\to A\ot_\ell A \\
&_\lambda T: A\ot_t A\to A\ot_r A
\qquad\quad\text{and}\qquad\quad
_\rho T: A\ot^s A\to A\ot_r A.
\endalign$$

\bf\snl Proof\rm:
Consider the canonical map $T_1:A\ot A\to A\ot A$, given by 
$$T_1(a\ot b)=\Delta(a)(1\ot b).$$
We know that the range of this map is $E(A\ot A)$ and that the kernel is \newline 
$(A\ot 1)(1-F_1)(1\ot A)$. On the other hand, from Proposition 2.8 we know that the projection map $\pi_\ell:A\ot A\to A\ot_\ell A$, restricted to $E(A\ot A)$ is an isomorphism. Similarly from Proposition 2.10 we know that the projection map $\pi_s:A\ot A\to A\ot_s A$, restricted to $(A\ot 1)(1-F_1)(1\ot A)$ is also an isomorphism. It follows that we have a bijective map $T_\rho:A\ot_s A\to A\ot_\ell A$ making the following diagram commute:
$$\CD
A\ot A @>T_1> > A\ot A\\
@V \pi_s V  V @V  V \pi_\ell V \\
A\ot_s A @>T_\rho> > A\ot_\ell A
\endCD$$
\snl
A similar proof can be given for the other cases, using the appropriate results obtained in Proposition 2.8 and Proposition 2.10.
\hfill$\square$
\einspr

Compare this result with the formulas found after Notation 6.3 in Section 6 of [T-VD].
\snl
This all together means that we have a left multiplier bialgebroid as well as a right multiplier bialgebroid with bijective canonical maps. As also the necessary compatibility conditions are satisfied, we get the following.

\iinspr{3.5} Theorem \rm 
The left and the right multiplier bialgebroids we have constructed combine to a regular multiplier Hopf algebroid in the sense of Definition 6.4 of [T-VD].
\hfill$\square$
\einspr

\it The counit and the antipode of the associated multiplier algebroid \rm
\nl
It follows from the general theory of regular multiplier Hopf algebroids that a counit exists; see Proposition 4.8 of [T-VD]. In this case we can give an explicit formula in terms of the data of the original regular weak multiplier Hopf algebra.
\snl
Recall that we use $\varepsilon_s$ and $\varepsilon_t$ for the source and target maps, from $A$ to $B$ and $C$ respectively, defined for the original weak multiplier Hopf algebra $(A,\Delta)$. For the left and right counits, as considered in the theory of weak multiplier Hopf algebroids, we will use here $\varepsilon_B$ and $\varepsilon_C$ respectively. 

\iinspr{3.6} Proposition \rm 
If $(A,\Delta)$ is a regular weak multiplier Hopf algebra, then the left counit $\varepsilon_B$ of the associated left multiplier bialgebroid $((B,A,\iota,S_B),\Delta_B)$ is the map from $A$ to $B$ defined by
$$\varepsilon_B(a)=\sum_{(a)}a_{(2)}S^{-1}(a_{(1)})$$
where we use the Sweedler notation for the original coproduct $\Delta$ and where $S$ is the original antipode of $(A,\Delta)$.

\bf\snl Proof\rm: 
Define $\varepsilon_B$ as above. We see that $\varepsilon_B(a)=S^{-1}(\varepsilon_t(a))$ and this belongs to $S^{-1}(C)$. It follows that $\varepsilon_B$ is a map from $A$ to $B$. 
We have 
$$(\varepsilon_B\ot\iota)(\Delta(a)(1\ot b))
=\sum_{(a)} a_{(2)}S^{-1}(a_{(1)}) \ot a_{(3)}b $$
in $A\ot A$ for all $a,b\in A$. Because $\varepsilon_B(xa)=x\varepsilon_B(a)$ for $x\in B$ and $a\in A$, it is easily seen that, in the algebroid framework, this formula reads as
$$(\varepsilon_B\ot\iota)T_\rho(a\ot b) 
=\sum_{(a)} a_{(2)}S^{-1}(a_{(1)}) \ot a_{(3)}b\tag"(3.1)"$$
where $T_\rho(a\ot b)=\Delta_B(a)(1\ot b)$ in $A\ot_\ell A$ (cf.\ Notations 3.3 above) and where the expression (3.1) is considered in the balanced tensor product $B\ot_\ell A$, defined via the relations $xx'\ot a=x'\ot S_B(x)a$ for $a\in A$ and $x,x'\in B$. This space is naturally identified with $CA\subseteq A$ via the map 
$$x\ot a\mapsto S_B(x)A.$$ 
If we apply this identification map to (3.1) above, we end up with
$$\sum_{(a)} a_{(1)}S(a_{(2)})a_{(3)}b=ab.$$
This shows that $\varepsilon_B$ satisfies the first diagram in Definition 4.1 of [T-VD]. So, it is the left counit for the left multiplier bialgebroid $((B,A,\iota,S_B),\Delta_B)$.
\hfill$\square$
\einspr

Remark that the conditions in Proposition 4.9 of [T-VD] are fulfilled. Indeed, the map $\varepsilon_B$ itself satisfies
$$\varepsilon_B(xa)=x\varepsilon_B(a)
\qquad\quad\text{and}\qquad\quad
\varepsilon_B(S_B(x)a)=\varepsilon_B(a)x$$
for all $x\in B$ and $a\in A$. We have defined $B$ as the set $\varepsilon_s(A)$ and we know that $S_B(B)=C=\varepsilon_t(A)$. Now because $\varepsilon_B(a)=S^{-1}(\varepsilon_t(a))$ we see that $B$ is also equal to the set $\varepsilon_B(A)$. It follows that the sets $I^s$ and $I^t$, as defined in Section 4 of [T-VD] are nothing else but $B$. And from Proposition 2.9 of [VD-W2], we know that $BA=CA=A$ and so the conditions in Proposition 4.8 of [T-VD] are indeed fulfilled.
\nl
We have a similar {\it right-handed} result.

\iinspr{3.7} Proposition \rm
The right counit $\varepsilon_C$ of the right multiplier bialgebroid \newline $((C,A,\iota,S_C),\Delta_C)$ is the map from $A$ to $C$ given by
$$\varepsilon_C(a)=\sum_{(a)}S^{-1}(a_{(2)})a_{(1)}.$$
\vskip-0.8cm \hfill$\square$
\einspr

Remark that we have here
$$\varepsilon_C(ay)=\varepsilon_C(a)y
\qquad\quad\text{and}\qquad\quad
\varepsilon_C(aS_C(y))=y\varepsilon_C(a)$$
for all $y\in C$ and $a\in A$. Moreover we have 
$$\align(\varepsilon_C\ot\iota)_\rho T(a\ot b)
&=(\varepsilon_C\ot\iota)((1\ot b)\Delta_C(a))\\
&=\sum_{(a)}S^{-1}(a_{(2)})a_{(1)}\ot ba_{(3)}\tag"(3.2)"
\endalign$$
in the balanced tensor product $C\ot^t A$, defined by the relations $yy'\ot a=y'\ot ay$ where $a\in A$ and $y,y'\in C$. If this space is identified with $A$ via $y\ot a\mapsto ay$, the expression (3.2) yields
$$\sum_{(a)}ba_{(3)}S^{-1}(a_{(2)})a_{(1)}=ba.$$
This is precisely proving the first diagram in Definition 5.12 of [T-VD].
\nl
\it What about the antipode? \rm 
It should be no surprise that the antipode of the multiplier Hopf algebroid is the same as the original antipode. And it is instructive to verify the various diagrams in Section 6 of [T-VD]. Consider e.g.\ the two diagrams in Definition 6.6 of [T-VD].
\snl
For the first one, we have 
$$\mu(S\ot \iota)T_\rho(a\ot b)=\sum_{(a)}S(a_{(1)})a_{(2)}b$$
where $S$ is the antipode of the algebroid and where we use the Sweedler notation for $\Delta_B(a)$. We use $\mu$ for the multiplication. On the other hand we have 
$$\mu(S_C\varepsilon_C\ot \iota)(a\ot b)=\sum_{(a)}S_C(S^{-1}(a_{(2)})a_{(1)})b=\sum_{(a)}S(a_{(1)})a_{(2)}b$$
where now $S$ is the original antipode and with the Sweedler notation for the original coproduct on $A$. We see that we get the same expressions.
\snl
For the second diagram, we have similar formulas. On the one hand, we have
$$\mu(\iota\ot S)_\lambda T(a\ot b)=\sum_{(b)}ab_{(1)}S(b_{(2)})$$
with the antipode of the algebroid and the Sweedler notation for $\Delta_C(b)$. On the other hand, we find
$$\mu(\iota\ot S_B\varepsilon_B)(a\ot b)=\sum_{(b)}aS_B(b_{(2)}S^{-1}(b_{(1)}))=\sum_{(b)}ab_{(1)}S(b_{(2)})$$
where now we have the original antipode, as well as the Sweedler notation for the original coproduct. Again we find the same expressions.
\nl
\nl

\bf 4. Multiplier Hopf algebroids arising from weak multiplier Hopf algebras \rm
\nl
In this section we will give {\it a necessary and sufficient condition} for a regular multiplier Hopf algebroid to come from a regular weak multiplier Hopf algebra via the procedure described in the previous section.
\snl
The starting point in this section is precisely as in the second item of Section 2 (Compatible quantum graphs). So $A$ is an algebra with local units. We have two commuting subalgebras $B$ and $C$, sitting nicely in $M(A)$ and we have anti-isomorphisms 
$$S_B:B\to C \qquad\quad \text{and} \qquad\quad S_C:C\to B.$$
 In other words, we have a compatible pair of a left and a right quantum graph 
$$\Cal A_B:=(B,A,\iota,S_B)
\qquad\quad\text{and}\quad\qquad
\Cal A_C:=(C,A,\iota,S_C)$$
as in Proposition 2.5 of Section 2. 
\snl
The basic extra assumption we will need for the main theorem (Theorem 4.11 below) is the following.

\inspr{4.1} Assumption \rm 
We require that the algebra $B$ above is separable Frobenius as in Definition 2.1 and that a separating linear functional $\varphi$ on $B$ can be chosen so that its modular automorphism coincides with the inverse of the automorphism $S_CS_B$ obtained as the composition of the given anti-isomorphisms $S_B:B\to C$ and $S_C:C\to B$.
\hfill$\square$\einspr

Before we continue, we need to add some remarks.

\inspr{4.2} Remark \rm
i) It will not be sufficient only to require that $B$ is separable Frobenius. It can occur that there is a separating linear functional $\varphi$ on $B$ as in Definition 2.1, but that it cannot be chosen so that its modular automorphism $\sigma$ is the inverse of the given automorphism $S_CS_B$ of $B$. Indeed, any automorphism $\sigma$ of $B$ can occur by making appropriate choices for $C$ and $S_B, S_C$. And if $\sigma$ does not leave the center invariant, it can not be the modular automorphism of a faithful functional. See also Example 5.3.ii) in the next section.
\snl
ii) The extra assumption only depends on the pair of the algebra $B$ and the given automorphism $S_CS_B$. It does not refer to how $B$ and $C$ are realized as commuting subalgebras of $M(A)$. 
\hfill$\square$\einspr

We now recall the result proven in Proposition 2.3. We also need to consider the remark made after the proof of this proposition. 
\snl
Suppose that $B$ satisfies the above assumption.  
Then there is a separability idempotent $E\in M(B\ot C)$ so that 
$S_B$ and $S_C$ are the antipodal maps, canonically associated with $E$. In other words, we have
$$E(x\ot 1)=E(1\ot S_B(x))
\qquad\quad\text{and}\qquad\quad
(1\ot y)E=(S_C(y)\ot 1)E$$
for all $x\in B$ and $y\in C$ where $S_B$ and $S_C$ are the given anti-isomorphisms from $B$ to $C$ and from $C$ to $B$ respectively.
\snl
Remark that also the converse is true (see Proposition 2.4). Therefore, we could as well have started with two non-degenerate algebras $B$ and $C$ and a separability idempotent $E\in M(B\ot C)$, together with the assumption that $B$ and $C$ are realized as commuting subalgebras in $M(A)$ as before.
\nl
\it The main result \rm
\nl
Now we assume that moreover $A$ carries two coproducts $\Delta_B$ and $\Delta_C$ so that the pair 
$$((\Cal A_B,\Delta_B), (\Cal A_C,\Delta_C))$$
is a regular multiplier Hopf algebroid. We will proceed from these assumptions and try to construct an underlying regular weak multiplier Hopf algebra $(A,\Delta)$.
\snl
The first step in the procedure is to associate two coproducts $\Delta$ and $\Delta'$ on $A$. 
The first one will be induced from $\Delta_B$ while the second one from $\Delta_C$. 
\snl
Recall that the left coproduct $\Delta_B$ is a map from $A$ to $A\overline\ot_\ell A$. So by definition we have elements
$$\Delta_B(a)(1\ot b)
\qquad\quad\text{and}\qquad\quad
\Delta_B(a)(c\ot 1)
$$
in $A\ot_\ell A$ for all $a$, $b$, $c$ in $A$. See e.g.\ the first item in Section 3. 
\snl
Now consider also $A\ot A$ with the two actions of $A$, given by right multiplication in the first and in the second factor. This gives rise to the extended module, denoted as $A\overline\ot A$, whose elements $\xi$ are characterized by the fact that 
$$\xi(1\ot b)
\qquad\quad\text{and}\qquad\quad
\xi(c\ot 1)$$
are elements in $A\ot A$ for all $b$ and $c$ in $A$. Here we can view $A\overline\ot A$ as sitting in $\text{End}(A\ot A)$ in the obvious way. This will clearly make $A\overline\ot A$ into an algebra. In fact, it is a subalgebra of the algebra $L(A\ot A)$ of left multipliers of $A\ot A$. The situation is similar as for $A\overline\ot_\ell A$, or rather $L_{\text{reg}}(A_B\ot A)$, as in Section 3.
\snl
In Proposition 2.8 we have seen that the map $\theta_\ell$ from $A\ot_\ell A$ to $A\ot A$, given by $\theta_\ell(a\ot b)=E(a\ot b)$ is well-defined and that $\pi_\ell\circ\theta_\ell$ is the identity map on $A\ot_\ell A$ where  $\pi_\ell$ is the quotient map from $A\ot A$ to $A\ot_\ell A$.
\snl
We are now ready to obtain the first coproduct $\Delta:A\to M(A\ot A)$.

\inspr{4.3} Lemma \rm
There is a homomorphism $\Delta:A\to A\overline\ot A$ given by
$$\Delta(a)(1\ot b)=\theta_\ell(\Delta_B(a)(1\ot b))
\quad\quad\text{and}\quad\quad
\Delta(a)(c\ot 1)=\theta_\ell(\Delta_B(a)(c\ot 1))\tag"(4.1)"$$
for $a,b,c\in A$. It satisfies
$\Delta(a)=E\Delta(a)=\Delta(a)E$
for all $a\in A$. 
It  is coassociative in the sense that
$$(\Delta\ot\iota)(\Delta(a)(1\ot b))(c\ot 1\ot 1)=
((\iota\ot\Delta)(\Delta(a)(c\ot 1))(1\ot 1\ot b)
\tag"(4.2)"$$
holds for all $a,b,c$ in $A$.

\snl\bf Proof\rm:
i) From the definition of the map $\theta_\ell$ it is clear that 
$$\theta_\ell(\xi(1\ot b))=\theta_\ell(\xi)(1\ot b)
\qquad\quad\text{and}\qquad\quad
\theta_\ell(\xi(c\ot 1))=\theta_\ell(\xi)(c\ot 1)$$
for all $\xi\in A\ot_\ell A$ and $b,c\in A$. On the other hand, because $\Delta_B(a)\in A\overline\ot_\ell A$, it follows that
$$(\Delta_B(a)(1\ot b))(c\ot 1)=(\Delta_B(a)(c\ot 1))(b\ot 1)$$
for all $a,b,c\in A$. If we now apply $\theta_\ell$ to this last equation, it will follow that $\Delta(a)$ is well-defined in $A\overline\ot A$ by the formulas in (4.1).
\snl
ii) By the definition of $\Delta(a)$ we have 
$\Delta(a)(1\ot b)=\theta_\ell(\Delta_B(a)(1\ot b))$
for all $a,b\in A$. Because the range of $\theta_\ell$ is precisely $E(A\ot A)$, it is clear that 
$E\Delta(a)(1\ot b)=\Delta(a)(1\ot b)$. This holds for all $b$ and so $E\Delta(a)=\Delta(a)$ for all $a$. On the other hand we have
$$\Delta_B(a)E(p\ot q)=\Delta_B(a)(p\ot S_B(E_{(1)})E_{(2)}q)=\Delta_B(a)(p\ot q)$$
for all $a,p,q\in A$ because $\Delta_B(a)\in L_{\text{reg}}(A_B\times A)$. As before, we are using the Sweedler type notation for $E$, the fact that the left leg of $E$ belongs to $B$ and that $S_B(E_{(1)})E_{(2)}=1$. See the remark at the end of the preliminary section. If we now apply $\theta_\ell$ to the above equation, we arrive at
$$\Delta(a)E(p\ot q)=\Delta(a)(p\ot q)$$
for all $a,p,q$ and hence $\Delta(a)E=\Delta(a)$ for all $a$.
\snl
iii) Now we show that $\Delta$ is a homomorphism. We take $a$, $a'$ and $b$ in $A$ and start with the formula
$$\Delta_B(aa')(1\ot b)=\Delta_B(a)\Delta_B(a')(1\ot b)\tag"(4.3)"$$
that holds in $A\ot_\ell A$. If we apply $\theta_\ell$ on the left hand side of this equation, we find $\Delta(aa')(1\ot b)$ by the definition of $\Delta$. We now also want to apply $\theta_\ell$ to the right hand side of (4.3). We claim that $\theta_\ell(\Delta_B(a)\xi)=\Delta(a)\theta_\ell(\xi)$ for all $a\in A$ and $\xi\in A\ot_\ell A$. Then we find by applying $\theta_\ell$ on (4.3) 
$$\align \Delta(aa')(1\ot b)
&=\theta_\ell (\Delta_B(a)\Delta_B(a')(1\ot b))\\
&=\Delta(a)\theta_\ell(\Delta_B(a')(1\ot b))\\
&=\Delta(a)\Delta(a')(1\ot b)
\endalign$$
and we are done. 
\snl
To prove the claim, observe that
$$\theta_\ell(\Delta_B(a)(p\ot q))=\Delta(a)(p\ot q)=\Delta(a)E(p\ot q)=\Delta(a)\theta_\ell(p\ot q)$$
for all $a\in A$ and $p,q\in A$. Remark that in the first and the last expression, we view $p\ot q$ in $A\ot_\ell A$ and in the second and the third, we consider it as sitting in $A\ot A$.
\snl
iv) Finally we prove coassociativity of $\Delta$ as in (4.2). As the first leg of $E$ belongs to the algebra $B$, we will have
$$(\Delta_B\ot\iota)(E(p\ot q))=(1\ot E)(\Delta_B\ot\iota)(p\ot q)$$
for all $p,q\in A$. Now take $a,b,c\in A$. Then
$$\align
(\Delta\ot\iota)(\Delta(a)(1\ot b))&(c\ot 1\ot 1)\\
&=(E\ot 1)(\Delta_B\ot\iota)(E(\Delta_B(a)(1\ot b))(c\ot 1\ot 1)\\
&=(E\ot 1)(1\ot E)(\Delta_B\ot\iota)(\Delta_B(a)(1\ot b))(c\ot 1\ot 1).
\endalign$$
Similarly we will find
$$\align
((\iota\ot\Delta)(\Delta(a)(c\ot 1))&(1\ot 1\ot b)\\
&=(1\ot E)(\iota\ot\Delta_B)(E(\Delta_B(a)(c\ot 1))(1\ot 1\ot b)\\
&=(1\ot E)(E\ot 1)(\iota\ot\Delta_B)(\Delta_B(a)(c\ot 1))(1\ot 1\ot b).
\endalign$$
As $(E\ot 1)(1\ot E)=(1\ot E)(E\ot 1)$ in $M(A\ot A\ot A)$, coassociativity of $\Delta$ follows from coassociativity of $\Delta_B$. 
\hfill $\square$ \einspr

The proof of coassociativity in item iv) above  needs more and finer arguments. We have e.g.\ made no distinction between $E(p\ot q)$ and $\theta_l(p\ot q)$ where in the first place, we consider $p\ot q$ as sitting in $A\ot A$ while in the second case it is considered as an element in $A\ot_\ell A$. On the other hand, the result is so obvious and to give a more precise argument will not really clarify the proof.
\snl
In a completely similar fashion, we can associate another coproduct $\Delta'$ on $A$. Now we use $A\overline\ot' A$ for the extended module of elements $\eta$ satisfying and characterized by the requirement that 
$$(1\ot b)\eta
\qquad\quad\text{and}\qquad\quad
(c\ot 1)\eta
$$
belong to $A\ot A$ for all $b,c\in A$. Recall that the right coproduct $\Delta_C$ is a map from $A$ to $A\overline\ot_r A$. We will now use the map $\theta_r:A\overline\ot_r A\to A\ot A$ given by $\theta_r(a\ot b)=(a\ot b)E$.

\inspr{4.4} Lemma \rm
There is a homomorphism $\Delta':A\to A\overline\ot' A$ given by
$$
(1\ot b)\Delta'(a)=\theta_r((1\ot b)\Delta_C(a))
\qquad\quad\text{and}\qquad\quad
(c\ot 1)\Delta'(a)=\theta_r((c\ot 1)\Delta_C(a))
$$
for $a,b,c\in A$. We have $\Delta'(a)=E\Delta'(a)=\Delta'(a)E$ for all $a$. The coproduct $\Delta'$ is coassociative in the sense that
$$(c\ot 1\ot 1)(\Delta'\ot\iota)((1\ot b)\Delta'(a))
=(1\ot 1\ot b)(\iota\ot\Delta')(c\ot 1)\Delta'(a))$$
holds for all $a,b,c$ in $A$.
\hfill $\square$ \einspr

The proof of Lemma 4.4 is completely similar as the one of Lemma 4.3. 
\nl
We would like that these two coproducts coincide. And because $\Delta(a)$ is a left multiplier and $\Delta'(a)$ a right multiplier of $A\ot A$, it then will follow that $\Delta$ maps $A$ into $M(A\ot A)$ for all $a\in A$. So we will eventually have that  $\Delta$ is a (regular) coproduct on $A$.
\snl
We will of course need the given relation between the left coproduct $\Delta_B$ and the right coproduct $\Delta_C$ of the multiplier Hopf algebroid. However, before we are able to use this result to obtain the equality of $\Delta$ and $\Delta'$ on $A$, we need the counit for each of these two coproducts. They are obtained in the Lemma 4.5 below. We will also need {\it an extra assumption}. This will be discussed later, see Proposition 4.10 and the subsequent remarks there.
\snl
First recall the counital maps $\varepsilon_B:A\to B$ and $\varepsilon_C:A\to C$ obtained in Theorem 6.8 of [T-VD] for the given multiplier Hopf algebroid. They satisfy, among other properties, the following equalities:
$$\align
&\varepsilon_B(xa)=x\varepsilon_B(a)
\quad\qquad\text{and}\qquad\quad
\varepsilon_B(S_B(x)a)=\varepsilon_B(a)x \\
&\varepsilon_C(ay)=\varepsilon_C(a)y
\quad\qquad\text{and}\qquad\quad
\varepsilon_C(aS_C(y))=y\varepsilon_C(a) 
\endalign$$
where $a\in A$, $x\in B$ and $y\in C$. See a remark before Definition 6.6 of [T-VD].
\snl
We will need this to obtain the counits in the following lemma.

\inspr{4.5} Lemma \rm
Define $\varepsilon:A\to \Bbb C$ by $\varepsilon=\varphi_B\circ \varepsilon_B$ and $\varepsilon':A\to \Bbb C$ by $\varepsilon'=\varphi_C\circ\varepsilon_C$ where $\varphi_B$ and $\varphi_C$ are the linear functionals on $B$ and $C$ respectively characterized by the formulas
$$(\varphi_B\ot\iota)E=1
\qquad\quad\text{and}\qquad\quad
(\iota\ot\varphi_C)E=1.
$$
Then $\varepsilon$ is a counit for $\Delta$ and $\varepsilon'$ is a counit for $\Delta'$.

\snl\bf Proof\rm:
i) First consider the characterizing property of the left counit $\varepsilon_B$ as in the first diagram of Definition 4.1 of [T-VD]. The map $x\ot q\mapsto S_B(x)q$ from $B\ot A$ to $A$ will map $(\varepsilon_B\ot\iota)T_\rho(a\ot b)$ to $ab$ for all $a,b\in A$. Recall that 
$T_\rho(a\ot b)=\Delta_B(a)(1\ot b)$. It follows that 
$$\align
(\varepsilon\ot\iota)(\Delta(a)(1\ot b))
&=(\varphi_B\ot\iota)(\varepsilon_B\ot\iota)(E(\Delta_B(a)(1\ot b)))\\ 
&=(\varphi_B\ot\iota)(E(\varepsilon_B\ot\iota)(\Delta_B(a)(1\ot b))).
\endalign$$ 
We have used that the left leg of $E$ belongs to $B$ and that $\varepsilon_B(xp)=x\varepsilon_B(p)$ when $x\in B$ and $p\in A$. If we write $(\varepsilon_B\ot\iota)(\Delta_B(a)(1\ot b))$ as $\sum_i x_i\ot q_j$, we see that 
$$
E(\varepsilon_B\ot\iota)(\Delta_B(a)(1\ot b))
=\sum_i E(x_i\ot q_i)
=\sum_i E(1\ot S_B(x_i)q_i)
$$
and we know that this is $E(1\ot ab)$. Therefore, if we apply $\varphi_B\ot\iota$ we get $ab$ and hence all together we find
$$(\varepsilon\ot\iota)(\Delta(a)(1\ot b))=ab$$
for all $a,b$ in $A$. 
\snl
ii) Next consider the second diagram in Definition 4.1 of [T-VD]. Now the map $p\ot x\mapsto xp$ from $A\ot B$ to $A$ will map $(\iota\ot\varepsilon_B)T_\lambda(c\ot a)$ to $ac$ for all $a,c\in A$. Recall that $T_\lambda(c\ot a)=\Delta_B(a)(c\ot 1)$. It follows that 
$$\align
(\iota\ot \varepsilon)(\Delta(a)(c\ot 1))
&=(\iota\ot\varphi_B)(\iota\ot\varepsilon_B)(E(\Delta_B(a)(c\ot 1)))\\
&=(\iota\ot\varphi_C)(\iota\ot S_B\circ\varepsilon_B)(E(\Delta_B(a)(c\ot 1)))
\endalign$$
Now $$S_B(\varepsilon_B(yp))=S_B(\varepsilon_B(p)S_B^{-1}(y))=yS_B(\varepsilon_B(p))$$
for all $y\in C$ and $p\in A$. Using now that the right leg of $E$ is in $C$, we find that the above expression equals
$$(\iota\ot\varphi_C)(E(\iota\ot S_B\circ\varepsilon_B)(\Delta_B(a)(c\ot 1))).$$
We now write $(\iota\ot\varepsilon_B)(\Delta_B(a)(c\ot 1))$ as $\sum_i p_i\ot x_i$. Then we see that 
$$E(\iota\ot S_B\circ\varepsilon_B)(\Delta_B(a)(c\ot 1))
=\sum_i E(p_i\ot S_B(x_i))=\sum E(x_ip_i \ot 1)$$
and this is $E(ac\ot 1)$. If we apply $\iota\ot\varphi_C$ we find $ac$ and all together we get
$$(\iota\ot \varepsilon)(\Delta(a)(c\ot 1)=ac$$
for all $a,c\in A$.
\snl
We see from i) and ii) that $\varepsilon$ is a counit for the coproduct $\Delta$ on $A$. 
\snl
A completely similar argument will give that $\varepsilon'$, defined as $\varphi_C\circ\varepsilon_C$, is a counit for the coproduct $\Delta'$ on $A$. Here we will need the module properties of $\varepsilon_C$ as recalled before the formulation of the lemma, as well as the first two diagrams in Definition 5.12 of [T-VD].
\hfill $\square$ \einspr

Again, we should be more precise as we indicated already after the proof of Lemma 4.3. Moreover we are not even working with genuine coproducts as we just have that $\Delta(a)$ is a left multiplier and $\Delta'(a)$ is a right multiplier of $A\ot A$ for any $a\in A$. Fortunately, we have taken care of these restrictions. We also do this for the next statement.
\snl
We know that the counits are unique provided the coproducts are full. Fullness of the coproducts $\Delta$ and $\Delta'$ will be a consequence of the following result.

\inspr{4.6} Lemma \rm
We have 
$$\align
\Delta(A)(1\ot A)&=\Delta(A)(A\ot 1)=E(A\ot A)\\
(1\ot A)\Delta'(A)&=(A\ot 1)\Delta'(A)=(A\ot A)E.
\endalign$$

\snl\bf Proof\rm:
Because we are working with a regular multiplier Hopf algebroid, all canonical maps are bijective by assumption (see Definition 6.4 of [T-VD]). In particular, the ranges of the maps $T_\rho$ and $T_\lambda$ are all of $A\ot_\ell A$. On the other hand, we know that the map $\theta_\ell$ is an isomorphism of $A\ot_\ell A$ to the subspace $E(A\ot A)$ of $A\ot A$. Combining these two results will give 
$$\Delta(A)(1\ot A)=E(A\ot A)
\quad\qquad\text{and}\qquad\quad
\Delta(A)(A\ot 1)=E(A\ot A).$$
Similarly, using that the maps $_\rho T$ and $_\lambda T$ have range all of $A\ot_r A$ and that $\theta_r$ is an isomorphism form $A\ot_r A$ to $(A\ot A)E$, we find
$$(1\ot A)\Delta'(A)=(A\ot A)E
\quad\qquad\text{and}\qquad\quad
(A\ot 1)\Delta'(A)=(A\ot A)E.
$$
\vskip -0.8cm
\hfill $\square$ \einspr

This indeed implies fullness of the coproducts. Assume e.g.\ that $V$ is a subspace of $A$ so that $\Delta(A)(1\ot A)\subseteq V\ot A$. Then we have $E(A\ot A)\subseteq V\ot A$. Because the left leg of $E$ is all of $B$, it follows that $BA\subseteq V$. But as we know that $BA=A$, we get $V=A$. Similarly for the right leg of $\Delta$ and for the legs of $\Delta'$. Hence these counits obtained in Lemma 4.5 are unique.
\snl
This property already takes care of the condition on the ranges of the canonical maps for a regular weak multiplier Hopf algebra. See Definition 1.14 in [VD-W1]. And observe that we want to obtain a {\it regular} weak multiplier Hopf algebra so that we have to consider not only the ranges of the canonical maps $T_1$ and $T_2$, but also of the maps $T_3$ and $T_4$. Of course, we still need to argue that $\Delta$ and $\Delta'$ coincide (see further).
\nl
Next we show that the coproducts have the expected behavior on the legs of $E$. This is condition ii) in Definition 1.14 of [VD-W1].

\inspr{4.7} Lemma \rm
We have 
$$(\iota\ot \Delta)(E)=(E\ot 1)(1\ot E)=(1\ot E)(E\ot 1)$$
and similarly for the other coproduct $\Delta'$.
\hfill $\square$ \einspr

Recall first that we have not yet obtained the equality of $\Delta$ and $\Delta'$. This implies not only that we have to prove these results for the two, but also that we have to be a bit more careful. At this moment, we only know that $\Delta(a)$ is a left multiplier and that $\Delta'(a)$ is right multiplier of $A\ot A$. This means that formulas should be given a meaning by multiplying, left or right, and in the appropriate way, with elements of $A$.
\snl
We know the behavior of $\Delta_B$ and $\Delta_C$ on $B$ and on $C$. See the formulas 6.1 in Section 6 of [T-VD]. This will imply the same for the coproducts $\Delta$ and $\Delta'$. We will e.g.\ get that $\Delta(y)=E(y\ot 1)$ for $y\in C$ and as $C$ is the right leg of $E$, this will imply
$$(\iota\ot\Delta)E=(1\ot E)(E\ot 1).$$
In fact, these two results are the same. All the other cases are obtained in a completely similar fashion. It is also possible to use that the two legs of $E$ commute. 
\nl
Next we look again at the canonical maps $T_1$, $T_2$, $T_3$ and $T_4$ associated with the coproducts $\Delta$ and $\Delta'$ on $A$ by the formulas
$$\align 
T_1(a\ot b) &=\Delta(a)(1\ot b)
\quad\quad\text{and}\quad\quad
T_2(a\ot b) = (a\ot 1)\Delta'(b)\\
T_3(a\ot b) &= (1\ot b)\Delta'(a)
\quad\quad\text{and}\quad\quad
T_4(a\ot b) = \Delta(b)(a\ot 1)
\endalign
$$
where $a,b\in A$. 
\snl
In Lemma 4.6 we have obtained the ranges of these maps. The range of $T_1$ and $T_4$ is $E(A\ot A)$ whereas the range of $T_2$ and  $T_3$ is $(A\ot A)E$. In the next lemma, we obtain the kernels of these maps.

\inspr{4.8} Lemma \rm We have 
$$\align\text{Ker}(T_1)&=(A\ot 1)(1-F_1)(1\ot A)
\quad\quad\text{and}\quad\quad
\text{Ker}(T_2)=(A\ot 1)(1-F_2)(1\ot A)\\
\text{Ker}(T_3)&=(1\ot A)(1-F_3)(A\ot 1)
\quad\quad\text{and}\quad\quad
\text{Ker}(T_4)=(1\ot A)(1-F_4)(A\ot 1)
\endalign$$
where $F_1$, $F_2$, $F_3$ and $F_4$ are as in Proposition 2.9.

\snl\bf Proof\rm:
Consider the map $T_1$. It is related with the map $T_\rho$ in the following commuting diagram (see the proof of Proposition 3.4): 
$$\CD
A\ot A @>T_1> > A\ot A\\
@V \pi_s V  V @V  V \pi_\ell V \\
A\ot_s A @>T_\rho> > A\ot_\ell A
\endCD$$
So if $\xi\in A\ot A$ and if $T_1(\xi)=0$, we have $T_\rho(\pi_s(\xi))=0$ and because $T_\rho$ is assumed to be injective, we find that also $\pi_s(\xi)=0$. On the other hand, if $\pi_s(\xi)=0$, we have $\pi_\ell(T_1(\xi))=0$ and because $\pi_\ell$ is injective on the range of $T_1$ we must have $T_1(\xi)=0$. So we see that the kernel of $T_1$ is the same as the kernel of $\pi_s$. From Proposition 2.10 we know that the kernel of $\pi_s$ is equal to $(A\ot 1)(1-F_1)(1\ot A)$ where $F_1=(\iota\ot S_C)E$. This proves the first formula of the lemma.
\snl
A similar argument works in the three other cases. For $T_2$ we need the map $_\lambda T$, in the case of $T_3$ we need $_\rho T$ and for $T_4$, the relevant map is $T_\lambda$. See Notation 3.3 and the subsequent remark. In all cases, we again need the results of Proposition 2.12.
\hfill $\square$ \einspr

So far we have not considered any relation between the left and the right structures. We now proceed and try to show that the two coproducts $\Delta$ and $\Delta'$ coincide. First we need the joint coassociativity rules.

\inspr{4.9} Lemma \rm
For all $a,b,c\in A$ we have
$$\hbox{ }(c\ot 1\ot 1)(\Delta'\ot\iota)(\Delta(a)(1\ot b))
=(\iota\ot\Delta)((c\ot 1)\Delta'(a))(1\ot 1\ot b)\tag"(4.4)"$$
\vskip -40pt
\inspr{}
$$(\Delta\ot\iota)((1\ot b)\Delta'(a))(c\ot 1\ot 1)
=(1\ot 1\ot b)(\iota\ot\Delta')(\Delta(a)(c\ot 1)).\tag"(4.5)"
$$

\snl\bf Proof\rm: We start from the joint coassociativity for the left and the right coproducts of the original multiplier Hopf algebroid. The first one says that
$$(c\ot 1\ot 1)(\Delta_C\ot\iota)(\Delta_B(a)(1\ot b))
=(\iota\ot\Delta_B)((c\ot 1)\Delta_C(a))(1\ot 1\ot b)
\tag"(4.6)"
$$
holds in $A\ot_r A\ot_\ell A$ for all $a,b,c\in A$. See Definition 6.4 in [T-VD]. Here the triple balanced tensor product $A\ot_r A\ot_\ell A$ is defined by the relations
$$\align
p\ot xq\ot r &= p\ot q\ot S_B(x)r \\
p\ot qy\ot r &= pS_C(y)\ot q\ot r
\endalign$$
where $p,q,r\in A$ and $x\in B$ and $y\in C$.
\snl
We now proceed as in the proof of the coassociativity of the coproducts $\Delta$ and $\Delta'$ on $A$ (item iv) in the proof of Lemma 4.3). For all $a,b,c\in A$ we have
$$\align
(c\ot 1\ot 1)(\Delta'\ot\iota)(\Delta(a)(1\ot b))
&=(c\ot 1\ot 1)(\Delta_C\ot\iota)(E\Delta_B(a)(1\ot b))(E\ot 1)\\
&=(c\ot 1\ot 1)(1\ot E)(\Delta_C\ot\iota)(\Delta_B(a)(1\ot b))(E\ot 1)\\
&=(1\ot E)(c\ot 1\ot 1)(\Delta_C\ot\iota)(\Delta_B(a)(1\ot b))(E\ot 1)
\endalign$$
where we have used that the first leg of $E$ belongs to $B$. 
On the other hand we have 
$$\align
(\iota\ot\Delta)((c\ot 1)\Delta'(a))(1\ot 1\ot b)
&=(1\ot E)(\iota\ot\Delta_B)((c\ot 1)\Delta_C(a)E)(1\ot 1\ot b)\\
&=(1\ot E)(\iota\ot\Delta_B)((c\ot 1)\Delta_C(a))(E\ot 1)(1\ot 1\ot b)\\
&=(1\ot E)(\iota\ot\Delta_B)((c\ot 1)\Delta_C(a))(1\ot 1\ot b)(E\ot 1)
\endalign$$
where now we used that the right leg of $E$ belongs to $C$. So coassociativity as in (4.4) follows now from coassociativity as in (4.6).
\snl
Similarly coassociativity as in (4.5) will follow from the second joint coassociativity rule 
$$(\Delta_B\ot\iota)((1\ot b)\Delta_C(a))(c\ot 1\ot 1)
=(1\ot 1\ot b)(\iota\ot\Delta_C)(\Delta_B(a)(c\ot 1))
$$
in Definition 6.4 of [T-VD].
\hfill $\square$ \einspr

The remark made after the proof of Lemma 3.3 is also relevant here.
\snl
Finally we are now able to show equality of $\Delta$ and $\Delta'$. However we need {\it an extra condition}. We will comment on it after the proof of the next proposition.

\iinspr{4.10} Proposition \rm Assume that the counits $\varepsilon$ and $\varepsilon'$ are equal. Then also the coproducts $\Delta$ and $\Delta'$ coincide in the sense that $p(\Delta(a)q)=(p\Delta'(a))q$ for all $a\in A$ and $p,q\in A\ot A$. In particular, $\Delta$ will map $A$ to $M(A\ot A)$.

\bf \snl Proof\rm:
Consider coassociativity as in (4.4):
$$(c\ot 1\ot 1)(\Delta'\ot\iota)(\Delta(a)(1\ot b))
=(\iota\ot\Delta)((c\ot 1)\Delta'(a))(1\ot 1\ot b).
$$
If we apply $\varepsilon'$ on the middle leg of the left hand side of this equality, we find 
$$(c\ot 1)\Delta(a)(1\ot b).$$
If on the other hand, we apply $\varepsilon$ on the middle leg of the right hand side of this equality, we find
$$(c\ot 1)\Delta'(a)(1\ot b).$$
So when $\varepsilon=\varepsilon'$ it follows that $\Delta=\Delta'$.
\hfill $\square$ \einspr

One has $S_C(\varepsilon_C(a))=\varepsilon_B(S(a))$ for all $a$. The formula as such is not found in the original paper [T-VD], but it can be obtained from the various results in Section 6 of [T-VD]. It follows that 
$$\varepsilon(S(a)=\varphi_B(\varepsilon_B(S(a)))=\varphi_B(S_C(\varepsilon_C(a)))=\varphi_C(\varepsilon_C(a))=\varepsilon'(a)$$
for all $a$. Therefore $\varepsilon=\varepsilon'$ if and only if $\varepsilon$ is invariant under the antipode $S$.
\snl
Remark that this condition involves the linear functionals $\varphi_B$ and $\varphi_C$ on the base algebras $B$ and $C$, but not only that. Indeed, the counits $\varepsilon$ and $\varepsilon'$ are given by $\varphi_B\circ \varepsilon_B$ and $\varphi_C\circ \varepsilon_C$. So, the condition $\varepsilon=\varepsilon'$ involves also the counital maps $\varepsilon_B$ and $\varepsilon_C$ of the given multiplier Hopf algebroid. A similar remark can be made when the condition is formulated as the invariance of $\varepsilon$ under the antipode $S$.
\snl
It is possible to give an example of a regular multiplier Hopf algebroid with a base algebra $B$ that is separable Frobenius, such that the faithful linear functional $\varphi_B$ has the right modular automorphism, but where the associated counit $\varepsilon$ is not invariant under the given antipode $S$. In fact, such an example already exists in the framework of Hopf algebroids. We will give such an example in the next section, see Proposition 5.4 and 5.5.
\snl
We now formulate the main result of this section. For completeness and the convenience of the reader, we include all the details. For the notations and concepts in the formulation, we refer to the preceding.

\iinspr{4.11} Theorem \rm
Assume that $A$ is an algebra with local units and that $B$ and $C$ are commuting subalgebras of $M(A)$ satisfying
$AB=BA=A$ and $AC=CA=A$. Assume that there is a separability idempotent $E\in M(B\ot C)$. Denote by  
$$S_B:B\to C 
\qquad\quad\text{and} \qquad\quad
S_C:C\to B$$ 
the anti-isomorphisms characterized by 
$$E(x\ot 1)=E(1\ot S_B(x))
\qquad\quad\text{and}\qquad\quad
(1\ot y)E=(S_C(y)\ot 1)E
$$ 
for $x\in B$ and $y\in C$. We also denote by $\varphi_B$ and $\varphi_C$ the unique linear functionals on $B$ and $C$ satisfying
$$(\varphi_B\ot\iota)E=1
\qquad\quad\text{and}\qquad\quad
(\iota\ot\varphi_C)E=1.$$
Finally we assume the existence of coproducts 
$$\Delta_B:A\to A\ot_\ell A
\qquad\quad\text{and} \qquad\quad
\Delta_C:A\to A\ot_r A$$ 
that turn the compatible pair
$$\Cal A_B:=(B,A,\iota,S_B)
\qquad\quad\text{and}\quad\qquad
\Cal A_C:=(C,A,\iota,S_C)$$
of a left and a right quantum graph into a regular multiplier Hopf algebroid. 
\snl
Then there exist coproducts $\Delta$ and $\Delta'$ on $A$ with values in the left and the right multiplier algebra of $A\ot A$ respectively, defined by 
$$\Delta(a)(1\ot b)=E(\Delta_B(a)(1\ot b))
\qquad\quad\text{and}\qquad\quad
(c\ot 1)\Delta'(a)=((c\ot 1)\Delta_C(a))E
$$
where $a,b,c\in A$. These two coproducts are full in the sense that the left and right legs of each of them equals $A$. 
\snl
The composition $\varphi_B\circ \varepsilon_B$ is a counit $\varepsilon$ on $(A,\Delta)$ and the composition $\varphi_C\circ\varepsilon_C$ is a counit $\varepsilon'$ for $(A,\Delta')$. These counits $\varepsilon$ and $\varepsilon'$ are the same if and only if one of them is invariant for the antipode $S$ of the original Hopf algebroid. In that case, the coproducts $\Delta$ and $\Delta'$ coincide and the pair $(A,\Delta)$ is a regular weak multiplier Hopf algebra. The canonical idempotent is $E$ and the antipode $S$ of the regular weak multiplier Hopf algebras is the same as the original antipode $S$ of the multiplier algebroid.

\snl\bf Proof\rm:
i) We have seen that $(B,A,\iota,S_B)$ is a left quantum graph, that $(C,A,\iota,S_C)$ is a right quantum graph and that they are compatible. See Proposition 2.5 in Section 2.
\snl
ii) In Proposition 2.8 we have shown that there are well-defined maps
$$\theta_\ell:A\ot_\ell A \to A\ot A
\qquad\quad\text{and}\qquad\quad
\theta_r:A\ot_r A \to A\ot A
$$
given by
$$\theta_\ell(p\ot q)=E(p\ot q)
\qquad\quad\text{and}\qquad\quad
\theta_r(p\ot q)=(p\ot q)E
$$
for $p,q\in A$. The map $\theta_\ell$ is a bijection from $A\ot_\ell A$ to $E(A\ot A)$ and $\theta_r$ is bijective from $A\ot_r A$ to $(A\ot A)E$. They allow to define $\Delta$ and $\Delta'$ on $A$ by
$$\align
\Delta(a)(1\ot b)&=\theta_\ell(\Delta_B(a)(1\ot b))\\
(c\ot 1)\Delta'(a)&=\theta_r((c\ot 1)\Delta_C(a))
\endalign$$
for $a,b,c\in A$. Remark that $\Delta(a)$ is a left multiplier of $A\ot A$ and $\Delta'(a)$ is a right multiplier of $A\ot A$ for all $a$. In Lemma 4.3 and Lemma 4.4 we have seen that $\Delta$ and $\Delta'$ are  coproducts. As a consequence of Lemma 4.6, we could show that these coproducts are full in the sense that their two legs are all of $A$. 
\snl
iii) In Lemma 4.5 we have seen that $\varepsilon$ and $\varepsilon'$, defined on $A$ as $\varphi_B\circ\varepsilon_B$ and $\varphi_C\circ\varepsilon_C$ are counits for $\Delta$ and $\Delta'$ respectively.
\snl
iv) In Proposition 4.10 we showed that the equality $\varepsilon=\varepsilon'$ yields that also $\Delta=\Delta'$. In fact, this last equality means that $\Delta$ and $\Delta'$ satisfy e.g.\ 
$$(c\ot 1)(\Delta(a)(1\ot b))=((c\ot 1)\Delta'(a))(1\ot b)$$ 
for all $a,b,c\in A$. Therefore, $\Delta(a)$ and $\Delta'(a)$ determine the same multiplier of $A\ot A$. Then $\Delta:A\to M(A\ot A)$ is a regular and full coproduct on $A$ as defined in Definitions 1.1 and 1.4 of [VD-W1]. The functional $\varepsilon$ is a counit for this coproduct $\Delta$ as defined in Definition 1.3 of [VD-W1]. It is unique because $\Delta$ is full. 
\snl
We have argued in a remark following the proof of Proposition 4.11 that the condition $\varepsilon=\varepsilon'$ is equivalent with the requirement that $\varepsilon$ is invariant under the antipode $S$ of the multiplier Hopf algebroid. The same is of course true with the invariance of $\varepsilon'$ under $S$. 
\snl
v) In Lemma 4.6 we have seen that the range of the canonical map $T_1$ is $E(A\ot A)$ of $T_2$ is $(A\ot A)E$. This takes care of condition i) of Definition 1.14 in [VD-W1] (definition of a weak multiplier Hopf algebra). We want to show that $(A,\Delta)$ is regular and so we also need condition i) of Definition 1.14 in [VD-W1] for the pair $(A^{\text{op}},\Delta)$. For this we need that also also the range of the canonical map $T_3$ is $(A\ot A)E$ and that of $T_4$ is $E(A\ot A)$. This is also part of Lemma 4.6.
\snl
vi) In Lemma 4.7 we found 
$$(\iota\ot \Delta)(E)=(E\ot 1)(1\ot E)=(1\ot E)(E\ot 1).$$
This is condition ii) in Definition 1.14 of [VD-W1].
\snl
vi) Finally, we need to see that the kernels of the four canonical maps $T_1$, $T_2$, $T_3$ and $T_4$ all satisfy the requirement as in Definition 1.14 of [VD-W1]. This is essentially shown in Lemma 4.2. The result there is given in terms of the idempotents $F_i$ in $M(A\ot A^{\text{op}})$ whereas in Definition 1.14 of [VD-W1], the condition is in terms of the idempotent maps $G_i$. However, e.g.\ in Proposition 4.5 of [VD-W1], we find the relation between these two.
\snl
So we see that all together, we get a regular weak multiplier Hopf algebra $(A,\Delta)$.
\snl
To show that the antipode of the weak multiplier Hopf algebra is the same as the antipode of the original multiplier Hopf algebroid, one can give various arguments. At the end of the previous section, we have seen that the antipode of the multiplier Hopf algebroid associated to a given weak multiplier Hopf algebra is the same as the original antipode. Then the result follows because the antipode of a multiplier Hopf algebroid is unique.
\hfill $\square$ \einspr

Of course, to make this last argument complete, one has to argue that the two procedures, the one described in Section 3 (from a weak multiplier Hopf algebra to a multiplier Hopf algebroid) and the one developed in this section (from a multiplier Hopf algebroid to a weak multiplier Hopf algebra), are inverses of each other. However, this should be immediately clear from the two constructions.
\nl\nl

\bf 5. Examples and special cases \rm
\nl
In this section we will give examples and discuss some special cases. The aim is to illustrate various aspects of our results in the previous sections. 
\snl
Let $H$ be a group acting from the left on a space $X$ and denote the action as $h\times x\mapsto h.x$ for $h\in H$ and $x\in X$. The set of triples
$$G=\{ (y,h,x) \mid x,y\in X, h\in H \text{ and } y=h.x \}$$
is a groupoid when the product of two elements $(y,h,u)$ and $(v,k,x)$ is defined for $u=v$ by
$$(y,h,u)(u,k,x)=(y,hk,x)$$
where $x,y,u,v\in X$ and $h,k\in H$. This example of a groupoid illustrates very well the fact that the structure is a combination of two other ones, say that of a space and that of a group. Similarly, in the case of a weak multiplier Hopf algebra, as well as in that of a multiplier Hopf algebroid, there a two distinguished aspects. There is the base algebra on the one hand, corresponding to the 'space' aspect, and the Hopf-like structure on the other hand, corresponding to the 'group' aspect.
\snl
It is also clear from what we see in this paper that the relation of a weak multiplier Hopf algebra and its associated multiplier Hopf algebroid is situated mainly on the level of the base algebra, i.e.\ the 'space' ingredient. For this reason, we start our examples with weak multiplier Hopf algebras and multiplier Hopf algebroids where the second aspect, the 'group' ingredient, is trivial.
\nl
\it Examples with a trivial Hopf-like structure \rm
\nl
Recall first the following example of a regular weak multiplier Hopf algebra from Section 3 in [VD-W2].

\inspr{5.1} Example \rm 
Let $B$ and $C$ be non-degenerate idempotent algebras and assume that $E$ is a separability idempotent in $M(B\ot C)$. Let $A$ be the algebra $C\ot B$. We identify $B$ and $C$ with their natural images in $M(A)$ and hence we write $yx$ for $y\ot x$ in $A$ where (also further) $x\in B$ and $y\in C$. We have that the images of $B$ and $C$ in $M(A)$ commute and so $xy=yx$. There is a full and regular coproduct $\Delta$ on $A$ defined by 
$$\Delta(yx)=(y\ot 1)E(1\ot x)$$
for $x\in B$ and $y\in C$. The pair $(A,\Delta)$ is a regular weak multiplier Hopf algebra. The source and target maps are given by
$$\varepsilon_s(yx)=S_C(y)x
\qquad\quad\text{and}\qquad\quad
\varepsilon_t(yx)=yS_B(x).$$
The antipode is given by $S(yx)=S_B(x)S_C(y)$. The counit is given by $\varepsilon(yx)=\varphi_C(yS_B(x))$ and also by $\varepsilon(yx)=\varphi_B(S_C(y)x)$. For details we refer to Proposition 3.3 of [VD-W2].
\hfill $\square$ \einspr

We now obtain the following result.

\inspr{5.2} Proposition \rm
Consider the regular weak multiplier Hopf algebra $(A,\Delta)$ from Example 5.1. and the associated regular multiplier Hopf algebroid as in Section 3 of this paper. The spaces $A\ot_\ell A$ and $A\ot_r A$ are naturally identified with $C\ot A$ and $A\ot B$ respectively via the maps
$$y'x'\ot yx\mapsto y'\ot S_B(x')yx
\qquad\quad\text{and}\qquad\quad
y'x'\ot yx\mapsto y'x'S_C(y)\ot x$$
whenever $x,x'\in B$ and $y,y'\in C$. The left and right coproducts $\Delta_B$ and $\Delta_C$ then both take the form
$$\Delta_B(yx)=y\ot x
\qquad\quad\text{and}\qquad\quad
\Delta_C(yx)=y\ot x
$$
for $x\in B$ and $y\in C$ where in the first case $y\ot x$ is considered in $C\ot M(A)$ and in the second case in $M(A)\ot B$.
The counital maps $\varepsilon_B$ and $\varepsilon_C$ of the associated multiplier Hopf algebroid are given by
$$\varepsilon_B(xy)=xS_B^{-1}(y)
\qquad\quad\text{and}\qquad\quad
\varepsilon_C(xy)=S_C^{-1}(x)y$$
for $x\in B$ and $y\in C$.
\hfill $\square$ \einspr

The proof is straightforward. 
\snl
We can verify the formulas we get in Proposition 3.6 for the relation between the source and target maps $\varepsilon_s$ and $\varepsilon_t$, as we have in Example 5.1, and the counital maps $\varepsilon_B$ and $\varepsilon_C$ as found in Proposition 3.6 and 3.7 in Section 3 of this paper. Indeed we get 
$$\align
\varepsilon_B(xy)&=S^{-1}(\varepsilon_t(yx))=S_B^{-1}(yS_B(x))=xS_B^{-1}(y)\\
\varepsilon_C(xy)&=S^{-1}(\varepsilon_s(yx))=S_C^{-1}(S_C(y)x)=S_C^{-1}(x)y
\endalign$$ 
for $x\in B$ and $y\in C$
and we see that the formulas in Example 5.1 and in Proposition 5.2 are in accordance with the result of Proposition 3.6 and 3.7.
\snl
The above example is considered already in Section 8 of [T-VD].
\nl
We now use the above example to consider some examples of multiplier Hopf algebroids that do not arise from a weak multiplier Hopf algebra. For this we first remark that any pair of non-degenerate algebras $B$ and $C$ with anti-isomorphisms $S_B:B\to C$ and $S_C:C\ot B$ give rise to a regular multiplier Hopf algebroid as in Proposition 5.2. See Section 8 of [T-VD].

\inspr{5.3} Example \rm
i) Take any non-degenerate algebra and idempotent algebra $B$. Let $C=B^{\text{op}}$, take for $S_B$ the identity map and for $S_C$ its inverse. This will give a multiplier Hopf algebroid. If $B$ is not separable Frobenius, it can not come from a weak multiplier Hopf algebra.
\snl
ii) Again take a non-degenerate and idempotent algebra $B$ and $C=B^{\text{op}}$. Take again for $S_B$ the identity map, but now let $S_C=\sigma^{-1}S_B^{-1}$ where $\sigma$ is a given automorphism of $B$. Again this will give a regular multiplier Hopf algebroid. If $B$ is not separable Frobenius, it can not come from a weak multiplier Hopf algebra.
\snl
iii) Consider the situation as in ii). If $B$ is separable Frobenius, but if the automorphism $\sigma$ does not leave the center of $B$ invariant, the regular multiplier Hopf algebroid still has no underlying weak multiplier Hopf algebra. Indeed, for any faithful linear functional $\varphi_B$ on $B$ with a modular automorphism $\sigma_B$, this modular automorphism $\sigma_B$ will leave the center of $B$ invariant and so it can not be the given automorphism $\sigma$.  
\snl
iv) Finally, if $B$ is separable Frobenius, we know that any faithful linear functional will be a separating linear function.
Therefore if in iii) the automorphism is a modular automorphism of some faithful linear functional, then the given multiplier Hopf algebroid has an underlying weak multiplier Hopf algebra as in this paper. Also conversely, for such a weak multiplier Hopf algebra constructed from a pair of algebras $B$ and $C$ with anti-isomorphism $S_B$ and $S_C$ as in Example 5.1, we know that $B$ is separable Frobenius and that there is a separating functional $\varphi$ with a modular automorphism that coincides with the inverse of the automorphism $S_CS_B$ of $B$
\hfill $\square$ \einspr 

It is still not clear if, in the case of a separable Frobenius algebra $B$, any automorphism leaving the center invariant will be the modular automorphism of a separating linear functional. If this is the case, the situation as in Example 5.3 would be completely understood.
\nl 
\it More multiplier Hopf algebroids without an underlying weak multiplier Hopf algebra \rm
\nl
With the next example, we show that there are cases were the  base algebra of a regular multiplier Hopf algebroid satisfies all the required conditions but where the two coproducts $\Delta$ and $\Delta'$ do not coincide. Hence also in this case, the multiplier Hopf algebroid will not result from a weak multiplier Hopf algebra.
\snl
We start with any regular weak multiplier Hopf algebra $(A,\Delta)$. We assume that $u,v$ are invertible elements in the multiplier algebra $M(B)$ of the base algebra $B$ satisfying 
$$E(vu\ot 1)E=E. \tag"(5.1)"$$
This condition is fulfilled if $u$ and $v$ are each others inverses. However, there are other cases. Indeed, consider $v'=S_B(v)$ and $u'=S_C^{-1}(u)$ so that
$$E(vu\ot 1)E=E(v\ot 1)(u\ot 1)E=E(1\ot v')(1\ot u')E=E(1\ot v'u')E.$$
We see that condition (5.1) is also fulfilled if $u'$ and $v'$ are each others inverses. If $S_B$ and $S_C$ are each others inverses, this happens only if $u$ and $v$ are inverses of each other, but otherwise, it can happen that $u'$ and $v'$ are inverses of each other while $u$ and $v$ are not.
\snl
In [VD3] the following is proven. We use the notations as above.

\inspr{5.4} Proposition \rm
Define $\Delta'(a)=(u\ot 1)\Delta(a)(v\ot 1)$ for all $a$ in $A$. Then $(A,\Delta')$ is again a regular weak multiplier Hopf algebra. The canonical idempotent $E'$ of $(A,\Delta')$ is $(u\ot 1)E(v\ot 1)$. The counit $\varepsilon'$ for the new pair $(A,\Delta')$ is given by $a\mapsto \varepsilon(u^{-1}av^{-1})$ or equivalently by $a\mapsto \varepsilon({u'}^{-1}a{v'}^{-1})$  where $\varepsilon$ is the counit for $\Delta$. The new antipode $S'$ is given by $S'(a)=uS(vav^{-1})u^{-1}$ where $S$ is the antipode of the original pair $(A,\Delta)$. For the new source and target maps $\varepsilon'_s$ and $\varepsilon'_t$, we have
$$\varepsilon'_s(a)=u\varepsilon_s(u^{-1}a)
\qquad\quad\text{and}\qquad\quad
\varepsilon'_t(a)=\varepsilon_t(a{v'}^{-1})v'$$
for all $a$. 
\hfill $\square$ \einspr

Observe that the base algebra $B$ is the same for the two weak multiplier Hopf algebras $(A,\Delta)$ and $(A,\Delta')$. The same is true for the base algebra $C$. However, the associated anti-isomorphisms $S'_B$ and $S'_C$ are not the same as $S_B$ and $S_C$. Instead they are given by
$$S'_B(x)=S_B(vxv^{-1})
\qquad\quad\text{and}\qquad\quad
S'_C(y)=uS_C(y)u^{-1}$$
for $x\in B$ and $y\in C$.
\snl
In [VD3] it is also shown that we have {\it mixed coassociativity}. This means that
$$(\Delta'\ot \iota)\Delta(a)=(\iota\ot\Delta)\Delta'(a)
\qquad\quad\text{and}\qquad\quad
(\Delta\ot \iota)\Delta'(a)=(\iota\ot\Delta')\Delta(a)$$
for all $a\in A$. These two relations make sense in the multiplier algebra $M(A\ot A\ot A)$, either by covering with the necessary elements of $A$ of by extending the maps involved.
\snl
Then the two coproducts are used to obtain a regular multiplier Hopf algebroid (see Theorem 2.3 in [VD3]). We use the terminology from [T-VD]. 

\inspr{5.5} Proposition \rm
Consider the left quantum graph $\Cal A_B:=(B, A, \iota, S_B)$ associated with the original pair $(A,\Delta)$ and the right quantum graph $\Cal A_C:=(C, A, \iota, S'_C)$ associated with the new pair $(A,\Delta')$. Define a left coproduct $\Delta_B$ on $A$ as in i) of Proposition 3.2, using the coproduct $\Delta$ and define a right coproduct $\Delta_C$ on $A$ as in ii) of Proposition 3.2, using the new coproduct $\Delta'$. Then $((\Cal A_B,\Delta_B),(\Cal A_C,\Delta_C))$ is a regular multiplier Hopf algebroid. 
\hfill $\square$ \einspr

The proof is straightforward and details are found in [VD3]. The antipode $S_r$ of this multiplier Hopf algebroid is given by $S_r(a)=uS(a)u^{-1}$ for $a\in A$ where $S$ is the original antipode of $(A,\Delta)$ (see Proposition 2.5 in [VD3]).
\snl
It is also clear that the regular multiplier Hopf algebroid obtained in this way will not come from a single weak multiplier Hopf algebra, except if the two coproducts $\Delta$ and $\Delta'$ are the same, that is if the elements $u$ and $v$ are central and each others inverses.
\nl\nl

\bf 6. Conclusions and further research\rm
\nl
In Section 3 of this paper, we have seen how any regular weak multiplier Hopf algebra gives rise, in a natural way, to a regular multiplier Hopf algebroid. On the other hand, in Section 4, we have obtained various criteria for a regular multiplier Hopf algebroid to arise in this way.
\snl
There are several conditions and of a different nature. In the first place, the base algebra $B$ has to be separable Frobenius. This however is not sufficient. Indeed, the two anti-isomorphisms $S_B$ and $S_C$ are given, as well as the automorphism $\sigma_B$ of $B$, defined as the inverse of the composition $S_CS_B$ of the two anti-isomorphisms. If now this automorphism does not leave the center invariant, it can not be the modular automorphism of a faithful functional on $B$. This means that, even if $B$ is separable Frobenius, the given regular multiplier Hopf algebroid will not have an underlying weak multiplier Hopf algebra. We have seen in Section 5 that this situation can occur.
\snl
Besides of these conditions on the base algebras $B$ and $C$ and the given anti-isomorphisms $S_B$ and $S_C$, there is still another condition, of a completely different nature. Even if all the necessary conditions on these objects are fulfilled, the left and right coproducts $\Delta_B$ and $\Delta_C$ may give rise to different coproducts $\Delta$ and $\Delta'$ on the algebra $A$. Also in this case, there is no underlying weak multiplier Hopf algebra. Again, we have seen in Section 5 that this can occur.
\snl
This takes us to two interesting questions. First, is it possible to formulate a condition on the given multiplier Hopf algebroid, not related in any way with the conditions on the base algebra, that guarantees the equality of these coproducts if they exist? The second question goes in the other direction. Is it possible to develop a more general theory of weak multiplier Hopf algebras where two different coproducts $\Delta$ and $\Delta'$ are allowed, related with each other only by the joint coassociativity rules? Would any such theory make sense? We see that there are examples.
\snl
In other words, one may wonder if either the axioms of a weak multiplier Hopf algebra are too restrictive or if the axioms of a multiplier Hopf algebroid are too general?
\nl
Finally, there are two other aspects that need to be considered. First, there is the theory of integrals, both in the case of regular weak multiplier Hopf algebras and for regular multiplier Hopf algebroids. The relation between the two should not give any surprises. Secondly, there is the non-regular case that has to be examined. Of course, before this can be done, it would be necessary first to develop the more general theory of non-regular multiplier Hopf algebroids. On the other hand, the obvious link between the non-regular weak multiplier Hopf algebras and the multiplier Hopf algebroids, might be the source of inspiration to do this.

\newpage

\bf References \rm
\nl
{[\bf B]} G.\ B\"ohm: {\it An alternative notion of Hopf algebroid}. In "Hopf algebras in noncommutative geometry and physics",
pp. 31-53, Lecture Notes in Pure and Appl. Math., 239, Dekker, New York, 2005. See also arXiv:math/0301169 
\snl
{[\bf B-N-S]} G.\ B\"ohm, F.\ Nill \& K.\ Szlach\'anyi: {\it Weak Hopf algebras I. Integral theory and C$^*$-structure}. J.\ Algebra 221 (1999), 385-438. 
\snl
{[\bf K-VD]} B.-J.\ Kahng \& A.\ Van Daele: {\it The Larson-Sweedler theorem for weak multiplier Hopf algebras}. See arXiv: 1406.0299 [math.RA]
\snl
{[\bf M]} W.\ Murray: {Nakayama automorphisms of Frobenius algebras} J.\ of Algebra {\bf 269} (2003), 599-609.
\snl
{[\bf VD1]} A.\ Van Daele: {\it Tools for working with multiplier Hopf algebras}. ASJE (The Arabian Journal for Science and Engineering) C - Theme-Issue 33 (2008), 505--528.  See also arXiv: 0806.2089 [math.QA]
\snl
{[\bf VD2]} A.\ Van Daele: {\it Separability idempotents and multiplier algebras}. See arXiv: \newline 1301.4398v1 [math.RA] 
\snl
{[\bf VD3]} A.\ Van Daele: {\it Modified weak multiplier Hopf algebras}. Preprint University of Leuven (2014). 
\snl
{[\bf VD-W0]} A.\ Van Daele \& S.\ Wang: {\it Weak multiplier Hopf algebras. Preliminaries and basic examples}. Operator algebras and quantum groups. Banach Center Publication 98 (2012), 367-415. See also arXiv: 1210.3954 [math.RA].
\snl
{[\bf VD-W1]} A.\ Van Daele \& S.\ Wang: {\it Weak multiplier Hopf algebras I. The main theory}. Journal f\"ur die reine und angewandte Mathematik (Crelles Journal), to appear.  See also arXiv: 1210.4395v1 [math.RA]. 
\snl
{[\bf VD-W2]} A.\ Van Daele \& S.\ Wang: {\it Weak multiplier Hopf algebras II. The source and target algebras}. See arXiv: 1403.7906v1 [math.RA]
\snl
{[\bf T-VD]} T.\ Timmermann \& A.\ Van Daele: {\it Regular multiplier Hopf algebroids. Basic theory and examples}. See arXiv: 1307.0769v3 [math.QA]


\end